\documentclass[a4paper,11pt]{article}

\usepackage[pdftex]{hyperref}
\usepackage{stmaryrd}
\usepackage[normalem]{ulem}
\usepackage{graphicx,xcolor}
\usepackage{amsfonts}
\usepackage{amsmath,amssymb,mathrsfs,amsthm}
\usepackage[all]{xy} 
\usepackage{color}
\usepackage{times}
\usepackage{enumerate,vmargin}
\usepackage{verbatim}
\usepackage{paralist}
\usepackage[labelfont=bf,size=small,font=it]{caption}
\usepackage[numbers,sort&compress]{natbib}
\usepackage{tikz}
\usetikzlibrary{shapes,arrows}

\setpapersize{A4}
\setmarginsrb{2.5cm}{2cm}{2.5cm}{2cm}{.4cm}{4mm}{0cm}{7.5mm}

\hypersetup{
    unicode      = false,     
    pdftoolbar   = true,      
    pdfmenubar   = true,      
    pdffitwindow = true,      
    pdfnewwindow = true,      
    colorlinks   = true,      
    linkcolor    = blue,      
    citecolor    = blue,      
    filecolor    = blue,      
    urlcolor     = blue       
}

\theoremstyle{plain}
\newtheorem{lem}{Lemma}[section]
\newtheorem{thm}[lem]{Theorem}
\newtheorem{prop}[lem]{Proposition}
\newtheorem{cor}[lem]{Corollary}

\theoremstyle{definition}

\theoremstyle{remark}

\newcommand{\cR}{\mathcal R}

\newcommand{\bP}{\mathbf P}
\newcommand{\bE}{\mathbf E}

\newcommand{\bX}{\mathrm X}

\newcommand{\bDelta}{\boldsymbol\Delta}
\newcommand{\DA}{\bDelta^{\downarrow}_{\alpha}}

\newcommand{\rY}{\mathrm Y}
\newcommand{\re}{\mathrm e}

\newcommand{\bp}{\mathbf p}

\newcommand{\llb}{\llbracket}
\newcommand{\rrb}{\rrbracket}

\newcommand{\dhau}{\operatorname{d_{\mathrm{H}}}}
\newcommand{\dpr}{\operatorname{d_{\mathrm{Pr}}}}
\newcommand{\dgh}{\operatorname{d_{\mathrm{GH}}}}
\newcommand{\dgp}{\operatorname{d_{\mathrm{GP}}}}

\newcommand{\dghp}{\operatorname{d_{GHP}}}

\newcommand{\gr}{\operatorname{gr}}
\newcommand{\dgr}{d_{\gr}}

\newcommand{\Lf}{\operatorname{Lf}}
\newcommand{\Br}{\operatorname{Br}}

\newcommand{\supp}{\operatorname{supp}}

\newcommand{\br}{\operatorname{br}}

\newcommand{\nr}{\mathbf N_{\mathrm{nr}}}

\newcommand{\Ver}{\operatorname{Vervaat}}

\newcommand{\cT}{\mathcal T}

\newcommand{\bth}{\boldsymbol{\theta}}
\newcommand{\bTheta}{\boldsymbol\Theta}

\newcommand{\R}{\mathbb R}
\newcommand{\N}{\mathbb N}

\newcommand{\Var}{\operatorname{Var}}

\newcommand{\eqd}{\stackrel{(d)}=}

\title{\textsc{Stable trees as mixings of inhomogeneous continuum random trees}}

\author{Minmin \textsc{Wang}
\thanks{Department of Mathematics, University of Sussex, 
Falmer campus, Brighton, BN1 9QH, England, United Kingdom.  
Email: minmin.wang@sussex.ac.uk }
}

\date{\today}

\begin{document}

\maketitle
\begin{abstract}
It has been claimed in Aldous, Miermont and Pitman \cite{AMP04} that all L\'evy trees are mixings of inhomogeneous continuum random trees. We give a rigorous proof of this claim in the case of a stable branching mechanism, relying on a new procedure for recovering the tree distance from the graphical spanning trees that works simultaneously for stable trees and inhomogeneous continuum random trees. 

\smallskip
%
%
\end{abstract}


\section{Introduction}

\subsection{Background}

Continuum random trees are random metric spaces that appear in the scaling limits of finite trees. The most iconic example is the Brownian continuum random tree, initially introduced by Aldous \cite{aldcrt1} as the scaling limit of the uniform labelled trees. An intimate connection to the Brownian motion was revealed in \cite{LeG93} and it was shown in \cite{aldcrt3} to be the universal scaling limits for the $n$-vertex Bienaym\'e trees where the underlying offspring distribution has a finite variance. 
Since then, various generalisations to the Brownian continuum random tree have been invented. For our purpose here, we will focus on the following two cases:
\begin{itemize}
\item
{\bf L\'evy trees}, introduced by Le Gall and Le Jan \cite{LGLJ98} and extensively studied in Duquesne and Le Gall \cite{DuLG02}. This class of continuum random trees arises naturally in the large-size limit of general Bienaym\'e trees, and provides a geometric representation for the genealogies of continuous-state branching processes. 

\item
{\bf Inhomogeneous continuum random tree}, which appeared in the study of general additive coalescence by Aldous and Pitman \cite{Al00}, as well as an inhomogeneous version of the birthday problem by Camarri and Pitman \cite{Pi00}. This class of continuum random trees is also expected in the scaling limits of random trees with fixed degree sequences.

\end{itemize} 
Besides the aforementioned connections with branching processes, coalescence, discrete random trees, models of continuum random trees are also linked to fragmentation processes, critical random graphs, random planar maps, etc. Their widespread presence can be partly explained by the fact that their 
nice probabilistic properties make them a powerful tool in the study of more complicated problems. 


Unlike the Brownian continuum random tree, whose branch points are always binary, both L\'evy trees and inhomogeneous continuum random trees possess ``hubs'', namely branch points of infinite degrees. As it turns out, this shared feature is far from a simple coincidence, but points to a deep connection between the two classes of continuum random trees. Indeed, from their study on the exploration processes of the inhomogeneous continuum random trees, Aldous, Miermont and Pitman \cite{AMP04} claimed that L\'evy trees can be obtained from  inhomogeneous continuum random trees by taking a suitable mixing of the latter. 

Roughly speaking, exploration process of a continuum random tree is the analogue of the depth-first walk (also called Lukasiewicz walk) for a discrete tree. For L\'evy trees, the role of exploration process is played by spectrally positive L\'evy processes, or an excursion of those L\'evy processes if one considers a single tree. From these L\'evy processes, one can extract the distance of the corresponding L\'evy tree via the so-called height process. 

For the inhomogeneous continuum random trees, their exploration processes, as weak limits of the depth-first walks of the corresponding discrete trees, have been identified in \cite{AMP04} to be the Vervaat transforms of extremal exchangeable processes. However, the height process, which has played a crucial role in the above encoding of L\'evy trees by L\'evy processes, is only known to exist in special cases where a Brownian component is present (\cite{AMP04}). The construction of a height process for general exchangeable processes remains an open question.  In the absence of this height process as useful middleman, we can not navigate easily from the extremal exchangeable processes to inhomogeneous continuum random trees;  
instead we often rely on the so-called Line-breaking Algorithm of Aldous and Pitman \cite{Al00} to access information about the inhomogeneous continuum random trees.

Nevertheless, it does not take a huge leap of faith to believe that somehow hidden in those extremal exchangeable processes are all the information we need to build an inhomogeneous continuum random tree, based on our experiences with  discrete trees and L\'evy trees. Moreover, there should also be a uniform way to define the height process that works for both L\'evy trees and inhomogeneous continuum random tree. Together with Kallenberg's Theorem \cite{Kallenberg73} on the characterisation of exchangeable processes on $[0, 1]$, this belief has led to the following paragraph in \cite{AMP04}:
{\it
\begin{quote}
Our work suggest that there are many similarities between ICRTs and L\'evy trees. In fact, L\'evy trees turn out to be ``mixings'' of ICRTs in an analogous way that L\'evy bridges are mixing of extremal bridges with exchangeable increments. This will be pursued elsewhere.
\end{quote}}
\noindent
To the author's best knowledge, there has not been a proof for the statement. The current work aims to partially confirm this statement in the case of stable trees, although we believe that some of the proof ideas are robust enough to extend to general L\'evy trees.  

The paper is organised as follows: after introducing the stable trees and inhomogeneous continuum random trees in Sections~\ref{sec:stable-intro} and \ref{sec: main}, we announce the main result in the latter half of Section~\ref{sec: main}. We then describe the main steps of the proof, whose details are found in the remaining sections. 

\subsection{Stable tree}
\label{sec:stable-intro}

We fix some $\alpha\in (1, 2)$. Let $\mathrm Y^{\alpha}=(Y^{\alpha}_{t})_{t\ge 0}$ be a spectrally positive $\alpha$-stable process that has the Laplace transform:
\[
\psi(\lambda):=\log(\mathbb E[e^{-\lambda Y^{\alpha}_{1}}]) =\lambda^{\alpha}, \quad \text{for } t, \lambda> 0.
\]
We denote by $\rY^{\mathrm{br}}=(Y^{\mathrm{br}}_{t})_{0\le t\le 1}$ the bridge process of $\rY^{\alpha}$ that ends on $0$ at time $1$ (see Section~\ref{sec: stable proc} for a definition). We then build an excursion-type process $\bX^{\alpha}$ from $\rY^{\mathrm{br}}$ using the Vervaat transformation: let $\rho=\inf\{t: Y^{\mathrm{br}}_{t-}=\inf_{s\in [0, 1]}Y^{\mathrm{br}}_{s}\}$ be the first infimum point of $\rY^{\mathrm{br}}$, and  set
\begin{equation}
\label{def: Vervaat}
X^{\alpha}_{t} = 
\left\{
\begin{array}{ll}
Y^{\mathrm{br}}_{t+\rho}-Y^{\mathrm{br}}_{\rho-}\,,  &  0 \le t \le 1-\rho,\\
Y^{\mathrm{br}}_{t+\rho-1}-Y^{\mathrm{br}}_{\rho-}\,, & 1-\rho\le t\le 1.
\end{array}
\right.
\end{equation}
Note that $X^{\alpha}_{t}\ge 0$ for all $t\in [0, 1]$. In fact, the process $\bX^{\alpha}$ has the same distribution as an excursion of $\rY^{\alpha}$ above its running infimum conditioned on returning to $0$ at time $1$ (\cite{Chau97}). In consequence, the excursion theory for L\'evy processes allows us to transfer various results on $\rY^{\alpha}$ to $\bX^{\alpha}$. A particularly important application for us is the construction of 
height processes by Le Gall \& Le Jan \cite{LGLJ98} (see also~\cite{DuLG02}). Note that $\mathrm Y^{\alpha}$ satisfies Grey's condition: $\int^{\infty} d\lambda/\psi(\lambda) <\infty$. 
It follows that there exists a {\it continuous} process $\mathrm H=\mathrm H(\bX^{\alpha})=(H_{t})_{0\le t\le 1}$ characterised as follows: for each $t\in [0, 1]$, we have 
\begin{equation}
\label{def: H}
H_{t} =\frac{\Gamma(2-\alpha)}{\alpha}\lim_{\epsilon\to 0}\epsilon^{\alpha-1}\,\#\Big\{s\in (0, t]: X^{\alpha}_{s-}< \inf_{u\in [s, t]}X^{\alpha}_{u}, \,\Delta X^{\alpha}_{s}\ge \epsilon\Big\},
\end{equation}
where the limit exists in probability. See Eq.~(4.5) in \cite{LGLJ98}. 
Regarding $\mathrm H$ as a curve depicting the ``contour'' of a tree, we then extract the $\alpha$-stable tree in the following way. For each pair $(s, t)\in [0, 1]^{2}$, we introduce a symmetric function
\begin{equation}
\label{def: H-encoding}
d_{\alpha}(s, t) = H_{s}+H_{t}-2m(\mathrm H, s, t), \quad \text{where} \quad m(\mathrm H, s, t) = \inf\{H_{u}: s\wedge t \le u\le s\vee t\}.
\end{equation}
It can be readily checked that $d_{\alpha}$ defines a pseudo-metric on $[0, 1]$. To turn this into a true metric, we say $s\sim t$ if and only if $d_{\alpha}(s, t)=0$. Then $d_{\alpha}$ induces a metric on the quotient space $\cT_{\alpha} := [0, 1]/\sim$, which we still denote as $d_{\alpha}$. The {\it $\alpha$-stable tree} is the (random) metric space
\[
(\cT_{\alpha}, d_{\alpha}).
\]
This is a ``tree-like'' metric space in the sense that every pair is joined by a unique path which turns out to be a geodesic, i.e.~a {\em real tree}. 
More precisely, we regard the stable tree as a random element taking values in the space $\mathbb T$ of measured real trees, which is a Polish space under the so-called Gromov--Prokhorov topology. We defer the formal introduction of real trees and the Gromov--Prokhorov topology to Section \ref{sec: GH}. 
It is often convenient to consider $\cT_{\alpha}$ as a rooted tree, with the {\it root} $r_{\alpha}$ taken as the point $p(0)$ of $\cT_{\alpha}$, where $p: [0, 1]\to [0, 1]/\sim$ stands for the canonical projection. 
In addition, the stable tree is naturally equipped with a probability measure $\mu_{\alpha}$, defined as the pushforward of the Lebesgue measure on $[0, 1]$ by $p$. We will refer to $\mu_{\alpha}$ as the {\it mass measure} of $\cT_{\alpha}$. 

By analogy to the graph theory, we can also introduce a notion of {\it node degrees} for $\cT_{\alpha}$. More precisely, for 
 $v\in \cT_{\alpha}$, we denote $\deg(v, \cT_{\alpha})$ to be the number of connected components of $\cT_{\alpha}\setminus\{v\}$. We then classify the points of $\cT_{\alpha}$ into three categories: $v\in \cT_{\alpha}$ is a {\it leaf} if it has degree $1$; a {\it branch point} if $\deg(v, \cT_{\alpha})\ge 3 $ or $\deg(v, \cT_{\alpha})=\infty$; the rest of the nodes are all of degree $2$. We denote respectively the sets of leaves and of branch points of $\cT_{\alpha}$ by $\Lf(\cT_{\alpha})$ and $\Br(\cT_{\alpha})$. It can be shown that $\Br(\cT_{\alpha})$ is countably infinite while $\Lf(\cT_{\alpha})$ has the continuum cardinality; moreover, both sets are dense in $\cT_{\alpha}$ and the mass measure $\mu_{\alpha}$ is supported on $\Lf(\cT_{\alpha})$. 

As it turns out, every branch point in $\cT_{\alpha}$ has infinite degrees. To further discern their ``infiniteness'', we define the {\it local time} of the branch point $b\in \Br(\cT_{\alpha})$ as follows. Let $\{\cT^{b}(i): i\in \N\}$ be the collection of the connected components of $\cT_{\alpha}\setminus\{b\}$. Put 
\[
\mathrm{ht}(\cT^{b}(i))=\sup_{v\in \cT^{b}(i)}d_{\alpha}(v, b)
\]
to be the height of $\cT^{b}(i)$. The {\it local time} of $b$ is the following limit in the a.s.~sense:
\[
\Delta^{\alpha}(b) = \lim_{\epsilon\to 0} \frac{1}{v(\epsilon)}\#\{i\in \N:\mathrm{ht}(\cT^{b}(i))\ge \epsilon\}, \quad \text{ with } v(\epsilon)=\big( (\alpha-1)\epsilon\big)^{-\frac{1}{\alpha-1}}.
\]
See \cite{DuLeG05}, Theorem 4.7. Note that the distribution of the local times $\{\Delta^{\alpha}(b): b\in \Br(\cT_{\alpha})\}$ is known, thanks to the following result from \cite{DuLeG05}, which states a one-to-one correspondence between $\Br(\cT_{\alpha})$ and the collection of the jumps of $\bX^{\alpha}$. More precisely, recall the canonical projection $p$ from $[0, 1]$ to $\cT_{\alpha}$; 
then a.s.~for each $s\in [0, 1]$ with $\Delta X^{\alpha}_{s}>0$, we have $b=p(s)\in \Br(\cT_{\alpha})$ and 
\begin{equation}
\label{eq: corres}
\Delta^{\alpha}(b) = \Delta X_{s}^{ \alpha}. 
\end{equation}
Conversely, for each $b\in \Br(\cT_{\alpha})$, there is a unique jump time $s$ of $\bX^{\alpha}$ such that $b=p(s)$. 
In consequence, standard properties of the stable process imply that $\Delta^{\alpha}(b)$ are all distinct and satisfy
\[
\sum_{b\in \Br(\cT_{\alpha})} \big(\Delta^{\alpha}(b)\big)^{2}<\infty \text{ and } \sum_{b\in \Br(\cT_{\alpha})}\Delta^{\alpha}(b)=\infty \text{ almost surely.}
\]
Therefore, we can rank $\{\Delta^{\alpha}(b): b\in \Br(\cT_{\alpha})\}$ in non increasing order. Let us denote by $\bDelta^{\downarrow}_{\alpha}=(\Delta_{i})_{i\in \N}$ this re-ordering, which is itself a (random) element of the following set 
\begin{equation}
\label{def: bTheta}
\boldsymbol\Theta := \left\{\bth = (\theta_{i})_{i\in \N}: \theta_{1}\ge \theta_{2}\ge \cdots \ge 0,\, \sum_{i\in\N}\theta_{i}^{2}<\infty \text{ and } \sum_{i\in \N}\theta_{i}=\infty\right\}.
\end{equation}
Let us observe that if $\bth=(\theta_{i})_{i\in \N}\in \boldsymbol\Theta$, then necessarily $\theta_{i}>0$ for all $i\in \N$. We equip $\boldsymbol\Theta$ with the $\ell^{2}$-norm and refer to it as the {\it parameter space}.

\subsection{Aldous--Camarri--Pitman's Line-breaking Algorithm and the main result}
\label{sec: main}
Let $\bth=(\theta_{i})_{i\ge 1}$ be a non random element of $\boldsymbol\Theta$. 
Let $\|\bth\|$ stand for the $\ell_{2}$-norm of $\bth$, i.e.~$\|\bth\|^{2}=\sum_{i\ge 1}\theta_{i}^{2}$. 
The following {\it line-breaking} construction of the ICRT is a trivial extension to the original version presented in \cite{AlPi99, Pi00}, where it is assumed that $\|\bth\|^{2}=1$. 

\smallskip
\noindent
{\bf Line-breaking Algorithm. } 
Given the data $\bth=(\theta_{i})_{i\in\N}\in \boldsymbol\Theta$, we sample a collection of independent Poisson processes. 
For each $i\in \N$,  
let $\xi_{i, 1}<\xi_{i, 2}<\cdots $ be the jumps of a Poisson process on $\R_{+}=[0, \infty)$ with intensity $\theta_{i}$ per unit length. 
In the terminology of \cite{AlPi99, Pi00}, the points $\{\xi_{i, j}: j\ge 2, i\ge 1\}$ are referred to as the  \emph{cutpoints}. 
The fact that $\sum_{i}\theta_{i}^{2}<\infty$ ensures there is only a finite number of cutpoints in any finite interval. 
It follows that cutpoints can be ranked in an increasing order: let us denote by $\eta_{1}<\eta_{2}<\cdots$ this ranking. 
We further assign the colour $i$ to the point $\eta_{k}$ if and only if $\eta_{k}\in\{\xi_{i, j}:  j\ge 2\}$. 
To build a tree, we use these ranked cutpoints to partition the half-line $[0, \infty)$ into line segments $[\eta_{k}, \eta_{k+1}]$, $k\ge 0$, with the understanding that $\eta_{0}=0$. We then assemble these line segments into a tree by gluing the line segment $[\eta_{k}, \eta_{k+1}]$ to $\xi_{i, 1}$ if $\eta_{k}$ has colour $i$. Since $\xi_{i, 1}$ is less than any cutpoints of colour $i$, one can be convinced that this gluing procedure is well-defined. 
%
Let $\cR_{1}$ be  the single branch $[0, \eta_{1}]$. For $k\ge 2$, let $\cR_{k}$ be the real tree obtained after gluing $[\eta_{k-1}, \eta_{k}]$ to $\cR_{k-1}$. 
We obtain in this way an increasing sequence of metric spaces $(\cR_{k})_{k\ge 1}$. 
Let $d$ be the distance on $\cR_{k}$ induced by the Euclidean metric of $\R_{+}$.  
We then define $(\cT, d)$ to be the completion of $(\cup_{k}\mathcal R_{k}, d)$, which turns out to be a real tree. 
The root of this tree is set at $r:=0$.  
On the other hand, the role of {\it mass measure} is played by the following a.s.~limit in the weak topology of $(\cT, d)$:
\[
\mu:=\lim_{k\to\infty}\frac{1}{k}\sum_{i=1}^{k}\delta_{\eta_{k}}.
\]
The existence of the above limit is a consequence of Aldous' theory on continuum random trees \cite{aldcrt1, aldcrt3}.  
We shall refer to the triplet $(\cT, d, \mu)$ as the $\bth$-ICRT and denote by $\bP^{\bth}$ its law (and by $\bE^{\bth}$ expectations with respect to this law). 
Let us also note that the previous construction implies the following scaling property in $\bth$: for any $c>0$, we have
\begin{equation}
\label{id: icrt-scaling}
(\cT, d, \mu) \text{ under } \bP^{c\bth} \overset{(d)}{=} \Big(\cT, \tfrac1c\, d, \mu\Big) \text{ under }\bP^{\bth}.
\end{equation}

Let us recall the sequence $\bDelta^{\downarrow}_{\alpha}$ formed by the ranked jumps in a normalised stable excursion $\bX^{\alpha}$. 
The main result of the paper is the following  
\begin{thm}
\label{thm: main}
For any measurable functional $F: \mathbb T\to \mathbb R_{+}$, we have
\[
\int_{\boldsymbol\Theta}\mathbb P\big(\bDelta^{\downarrow}_{\alpha}\in d\bth\big)\,\bE^{\bth}\big[F(\cT, d, \mu)\big] = \mathbb E\big[F(\cT_{\alpha}, d_{\alpha}, \mu_{\alpha})\big]. 
\] 
\end{thm}

As an immediate consequence of Theorem \ref{thm: main}, we obtain a new construction for the stable tree.  
\begin{cor}
If we run the Line-breaking Algorithm with the data $\bDelta^{\downarrow}_{\alpha}$, then the continuum random tree obtained has the same distribution as $(\cT_{\alpha}, d_{\alpha}, \mu_{\alpha})$. 
\end{cor}


\medskip

\noindent
{\bf Notation.} Throughout the paper, we use the following uniform notation for graph trees and real trees. 
If $t$ is a (graph) tree and $v$ a vertex of $t$ (resp.~if $t$ is a real tree and $v\in t$), we denote by $\deg(v, t)$ the degree of $v$ in $t$. We also make the convention that  $\deg(v, t)=-\infty$ if $v\notin t$.  If $v, v'$ are two vertices of $t$, we denote by $\llb v, v'\rrb_{t}$ the unique path of $t$ connecting $v$ and $v'$.

\subsection{Outline of the proof}
\label{sec: outline}

\begin{figure}
\centering
\includegraphics[height = 3cm]{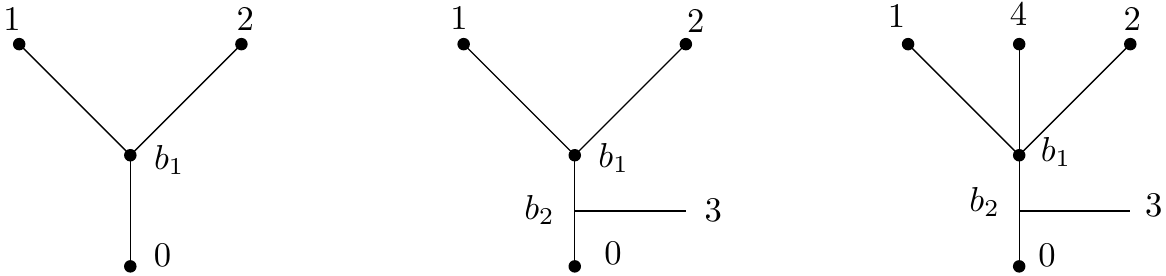}
\caption{\label{fig}An example of $T_{2}, T_{3}$ and $T_{4}$, together with the vertex labelling.}
\end{figure}

Our approach to the proof of Theorem~\ref{thm: main} is based upon a sequence of discrete approximations of the continuum random trees which works simultaneously for $\cT_{\alpha}$ and $\cT$. We first explain how this works for the stable tree. 
Let $(\sigma_{i})_{i\ge 1}$ be a sequence of i.i.d.~points of $\cT_{\alpha}$ with common distribution $\mu_{\alpha}$ and let $\cT_{k}$ be the smallest subtree of $\cT_{\alpha}$ that contains the root $r_{\alpha}$ and the first $k$ entries of $(\sigma_{i})_{i\ge 1}$. 
Note that $\cT_{k}\nearrow\cT$ as $k\to\infty$. 
On the other hand, we observe that $\cT_{k}$ has the ``shape'' of a discrete rooted tree with $k$ leaves, and we denote by $T_{k}$ this discrete tree (see Section~\ref{sec: stable-recov} for a more precise definition of $T_{k}$). 
Moreover, we regard $T_{k}$ as a labelled tree: the $k$ leaves, corresponding to the points $\sigma_{1}, \sigma_{2}, \dots, \sigma_{k}$ in $\cT_{k}$, are labelled from 1 to $k$; the branch points are labelled as $b_{1}, b_{2}, \cdots$, according to the order of their appearance in the sequence $(T_{k})_{k\ge 1}$; the root is relabelled as 0; see Fig.~\ref{fig} for an example. 
It is not difficult to see that this vertex labelling is consistence across $(T_{k})_{k\ge 1}$, so that $T_{k}$ appears as a subgraph of $T_{m}$, $m>k$.
In particular, this means that the vertex set $V(T_{k})$ of $T_{k}$ is a subset of $V(T_{k+1})$. Let us denote 
\[
V^{\alpha}_{\infty}=\bigcup_{k\ge 1}V(T_{k}).
\]

We claim that the $\alpha$-stable tree $(\cT_{\alpha}, d_{\alpha}, \mu_{\alpha})$ can be recovered from the sequence $(T_{k})_{k\ge 1}$ in a three-step procedure: first identify the local times of the branch points from the degree sequences of $(T_{k})_{k\ge 1}$, then recover the tree-distance $d_{\alpha}$ from these local times, and subsequently the mass measure $\mu_{\alpha}$. This is summarised in the following proposition, shown in Section \ref{sec: stable-recov}. 

\begin{prop}[Recovery of stable trees]
\label{prop: stable-recov}
The following statements hold true.
\begin{enumerate}[(i)]
\item
For each $v\in V^{\alpha}_{\infty}$\,, the following limit exists a.s.
\begin{equation}
\label{def: loc-stable'}
\widetilde\Delta^{\alpha}(v):=\lim_{k\to\infty} \frac{\deg(v, T_{k})}{k^{1/\alpha}}.
\end{equation}
Moreover, $\widetilde\Delta^{\alpha}(v)>0$ if and only if $v$ is a branch point of $T_{k}$ for some $k\ge 1$, and we have
\[
\big\{\widetilde\Delta^{\alpha}(v)>0: v\in V^{\alpha}_{\infty}\big\} = \{\Delta X^{\alpha}_{s}>0: 0\le s\le 1\} \quad\text{a.s.}
\] 
\item
For $v, v'\in V^{\alpha}_{\infty}$\,, let 
\[
N^{\alpha}_{\epsilon}(v, v')=\#\big\{w\in V^{\alpha}_{\infty}: \widetilde\Delta^{\alpha}(w)>\epsilon \text{ and } \exists\, k\in \N\text{ s.t. } w\in \llb v, v'\rrb_{T_{k}}\big\},
\]
which is a.s.~finite. The following limit exists in probability for all $v, v'\in V^{\alpha}_{\infty}$\,: 
\begin{equation}
\label{def: dist-stable}
\widetilde d_{\alpha}(v, v') =\frac{\Gamma(2-\alpha)}{\alpha}\cdot \lim_{\epsilon\to 0+} \epsilon^{\alpha-1}N^{\alpha}_{\epsilon}(v, v').
\end{equation}
Moreover, $\widetilde d_{\alpha}$ defines a metric on $V^{\alpha}_{\infty}$. Denoting $(\widetilde \cT_{\alpha}, \widetilde d_{\alpha})$ for  the completion of  $(V^{\alpha}_{\infty}, \widetilde d_{\alpha})$, we have $(\widetilde \cT_{\alpha}, \widetilde d_{\alpha})$ isometric to $(\cT_{\alpha}, d_{\alpha})$. 
\item
Let $\nu^{\alpha}_{k}$ be the uniform probability measure on the leaf set of $T_{k}$. Then $(\nu^{\alpha}_{k})_{k\ge 1}$ converges a.s.~to a limit $\widetilde \mu_{\alpha}$ and we have  $(\widetilde \cT_{\alpha}, \widetilde d_{\alpha}, \widetilde \mu_{\alpha})$ isometric to $(\cT_{\alpha}, d_{\alpha}, \mu_{\alpha})$. 
\end{enumerate}
\end{prop}

On the ICRT side, we similarly sample a sequence of subtrees of $\cT$: let $(\eta'_{i})_{i\ge 1}$ be a sequence of i.i.d.~points of $\cT$ with common law $\mu$; let $\cR'_{k}$ be the smallest subtree of $\cT$ containing the root $r$ and $\eta'_{1}$, $\eta'_{2}$, $\dots$, $\eta'_{k}$. In fact, as pointed out in \cite{Al00}, $(\cR'_{k})_{k\ge 1}$ has the same distribution as $(\cR_{k})_{k\ge 1}$, the sequence of real trees that appear in the Line-breaking Algorithm. 
Denote by $R'_{k}$ the ``shape'' of $\cR'_{k}$, which is a graph tree with no vertex of degree $2$. 
Denote by $V(R'_{k})$ the vertex set of $R'_{k}$ and by 
\[
V^{\bth}_{\infty}=\bigcup_{k\ge 1}V(R'_{k}).
\] 
The following result from Section \ref{sec: icrt-recov} says that the $\bth$-ICRT $(\cT, d, \mu)$ can be recovered from $(R'_{k})_{k\ge 1}$ in an analogous three-step procedure. 

\begin{prop}[Recovery of ICRTs]
\label{prop: icrt-recov}
The following statements hold true.
\begin{enumerate}[(i)]
\item
For each $v\in V^{\bth}_{\infty}$\,, the following limit exists in probability
\begin{equation}
\label{def: loc-icrt}
\Delta^{\bth}(v):=\lim_{k\to\infty}\frac{\deg(v, R'_{k})}{\Psi^{-1}_{\bth}(k)},  
\end{equation}
where $\Psi^{-1}_{\bth}$ is the inverse function of $\Psi_{\bth}(t)=\sum_{i\ge 1}(e^{-\theta_{i}t}-1+\theta_{i}t)$, the latter being strictly increasing as $\bth\in \bTheta$.
Moreover, $\Delta^{\bth}(v)>0$ if and only if $v$ is a branch point of $R'_{k}$ for some $k\ge 1$, and 
we have 
\[
\big\{\Delta^{\bth}(v)>0: v\in V^{\bth}_{\infty}\big\}=\{\theta_{i}: i\in \N\} \quad \text{a.s.}
\]
\item
For $v, v'\in V^{\bth}_{\infty}$\,, let 
\[
N^{\bth}_{\epsilon}(v, v')=\#\big\{w\in V^{\bth}_{\infty}: \Delta^{\bth}(w)>\epsilon \text{ and } \exists\, k\in \N\text{ s.t. } w\in \llb v, v'\rrb_{R'_{k}}\big\},
\]
which is finite since $\sum_{i}\theta_{i}^{2}<\infty$. 
The following limit exists in probability for all $v, v'\in V^{\bth}_{\infty}$\,: 
 \begin{equation}
\label{def: dist-icrt}
\widetilde d(v, v')=\lim_{\epsilon\to 0+}\frac{N^{\bth}_{\epsilon}(v, v')}{\gamma_{\bth}(\epsilon)} \,,
\end{equation}
where $\gamma_{\bth}(\epsilon)=\sum_{i}\theta_{i}\mathbf 1_{\{\theta_{i}>\epsilon\}}\nearrow \infty$ as $\epsilon\to 0$. 
Moreover, $\widetilde d$ defines a metric on $V^{\bth}_{\infty}$. Denoting $(\widetilde \cT, \widetilde d\,)$ for  the completion of  $(V^{\bth}_{\infty}, \widetilde d)$, we have $(\widetilde \cT, \widetilde d\,)$ isometric to $(\cT, d)$ under $\bP^{\bth}$. 
\item
Let $\nu^{\bth}_{k}$ be the uniform probability measure on the leaf set of $R'_{k}$. Then $(\nu^{\bth}_{k})_{k\ge 1}$ converges a.s.~to a limit $\widetilde \mu$ and we have  $(\widetilde \cT, \widetilde d, \widetilde \mu)$ isometric to $(\cT, d, \mu)$ under $\bP^{\bth}$. 
\end{enumerate}
\end{prop}

Let us point out the normalisations in \eqref{def: loc-stable'} and \eqref{def: dist-stable} are consistent with \eqref{def: loc-icrt} and \eqref{def: dist-icrt}. More precisely, we have the following result from Section \ref{sec: jumps}. 
\begin{lem}
\label{lem: normalisation}
Let $\bDelta^{\downarrow}_{\alpha}=(\Delta_{i})_{i\in \N}$ be the non increasing rearrangement of $\{\Delta X^{\alpha}_{s}>0: 0\le s\le 1\}$. Let 
\[
\Psi_{\bDelta^{\downarrow}_{\alpha}}(t)=\sum_{i\ge 1}(e^{-\Delta_{i}t}-1+\Delta_{i}t) \quad \text{and} \quad \gamma_{\bDelta^{\downarrow}_{\alpha}}(t)=\sum_{i\ge 1}\Delta_{i}\mathbf 1_{\{\Delta_{i}>t\}}, \quad t> 0. 
\]
Let $\Psi^{-1}_{\bDelta^{\downarrow}_{\alpha}}$ be the inverse function of $\Psi_{\bDelta^{\downarrow}_{\alpha}}$. Then we have the following limits in probability: 
\begin{equation}
\label{eq: asymp}
\frac{\Psi_{\bDelta^{\downarrow}_{\alpha}}(k)}{k^{\alpha}} \xrightarrow{k\to\infty} 1, \quad  \frac{\Psi^{-1}_{\bDelta^{\downarrow}_{\alpha}}(k)}{k^{1/\alpha}} \xrightarrow{k\to\infty} 1 \ \  \text{ and } \ \ 
\epsilon^{\alpha-1}\gamma_{\bDelta^{\downarrow}_{\alpha}}(\epsilon)\xrightarrow{\epsilon\to 0+} \frac{\alpha}{\Gamma(2-\alpha)}.  
\end{equation}
\end{lem}

We denote by $\mathbf T_{\mathrm{discrete}}$ the space of finite labelled (graph) trees equipped with the discrete topology and by $\mathbf T_{\mathrm{discrete}}^{\infty}$ the sequence of finite labelled trees equipped with the product topology.  
Lemma~\ref{lem: normalisation} combined with Propositions~\ref{prop: stable-recov} and~\ref{prop: icrt-recov} implies that we can find a {\it common} measurable function $\mathscr S: \boldsymbol\Theta\times\mathbf T_{\mathrm{discret}}^{\infty}\to \mathbb T$ so that we can write 
\[
(\cT_{\alpha}, d_{\alpha}, \mu_{\alpha})=\mathscr S\big(\boldsymbol\Delta^{\downarrow}_{\alpha}, (T_{k})_{k\ge 1}\big) \ \text{ and } \ (\cT, d, \mu)=\mathscr S\big(\bth, (R'_{k})_{k\ge 1}\big) \text{ under } \bP^{\bth}.
\]
Theorem \ref{thm: main} will then follow once we prove the next result. 
\begin{thm}
\label{thm: main'}
For any measurable functional $G: \mathbf T_{\mathrm{discret}}^{\infty}\to\R_{+}$, we have 
\[
\int_{\boldsymbol\Theta}\mathbb P\big(\bDelta^{\downarrow}_{\alpha}\in d\bth\big)\bE^{\bth}\big[G\big((R'_{k})_{k\ge 1}\big)\big] = \mathbb E\big[G\big((T_{k})_{k\ge 1}\big)\big]. 
\] 
\end{thm}

The proof of Theorem \ref{thm: main'}, given in Section~\ref{sec: main-pf}, relies upon the encoding of the stable tree and ICRT by certain stochastic processes with exchangeable increments. In the case of the stable tree, we have already encountered this coding process, which is 
the normalised excursion process $\mathrm X^{\alpha}$ of the stable process $\mathrm Y^{\alpha}$. 
Since the height process $\mathrm H$ is itself a functional of $\mathrm X^{\alpha}$, it should not come as a surprise that the spanning trees $(T_{k})_{k\ge 1}$ can be written as a measurable function of $\mathrm X^{\alpha}$ together with a sequence of i.i.d.~uniform points $(U_{i})_{i\ge 1}$ in $[0, 1]$. We give the explicit form of this function $\mathscr T^{\mathrm{lab}}$ in Section \ref{sec: coding} and show that 
\begin{equation}
\label{eq: sp-stable}
\{T_{k}: k\ge 1\} \eqd \big\{\mathscr T^{\mathrm{lab}}(\mathrm X^{\alpha}; \{U_{1}, \dots, U_{k}\}):  k\ge 1\big\}. 
\end{equation}
Let us note the encoding here can be seen as a ``coarser'' version of \eqref{def: H-encoding}, which only retains the shape of the trees but ignores the distances within. 

The candidate for the coding process of the $\bth$-ICRT has already been identified by Aldous, Miermont and Pitman \cite{AMP04}. It is closely connected to an extremal exchangeable  process $\mathrm Y^{\bth}=(Y^{\bth}_{t})_{0\le t\le 1}$ defined as 
\begin{equation}
\label{def: Ytheta}
Y^{\bth}_{t}= \sum_{i\ge 1}\theta_{i}(\mathbf 1_{\{\chi_{i}\le t\}}-t), \quad 0\le t\le 1,
\end{equation}
where in above $(\chi_{i})_{i\ge 1}$ is a sequence of i.i.d.~uniform points on $[0, 1]$, and the series on the right-hand side converges uniformly on $[0, 1]$ a.s. Using the Vervaat transform in \eqref{def: Vervaat} but replacing $\rY^{\mathrm{br}}$ with $\mathrm Y^{\bth}$, we can extract an excursion-type process $\bX^{\bth}$ and show in Section~\ref{sec: spnn-icrt} that 
\begin{equation}
\label{eq: sp-icrt}
\{R'_{k}: k\ge 1\} \eqd \big\{\mathscr T^{\mathrm{lab}}(\mathrm X^{\bth}; \{U_{1}, \dots, U_{k}\}):  k\ge 1\big\}. 
\end{equation}
Let us point out the proof of \eqref{eq: sp-icrt} is however quite different to that of \eqref{eq: sp-stable}. This is due to the lack of a height process for $\mathrm X^{\bth}$ which can play the same role as $\mathrm H$ for $\mathrm X^{\alpha}$. It then seems difficult to prove \eqref{eq: sp-icrt} directly based on the results in \cite{AMP04}. Instead, we go back to the discrete model (i.e.~$\bp$-trees) and introduce a discrete counterpart of \eqref{eq: sp-icrt}. We then work our way back through weak convergence arguments.

The final ingredient in the proof of Theorem~\ref{thm: main'} is provided by a celebrated theorem of Kallenberg \cite{Kallenberg73}, which implies that $\mathrm X^{\alpha}$ conditional on its jumps $\DA=\bth$ is distributed as $\mathrm X^{\bth}$. Together with \eqref{eq: sp-stable} and \eqref{eq: sp-icrt}, this will complete the proof of Theorem~\ref{thm: main'}, as we shall see in Section~\ref{sec: main-pf}. 

%

\section{Some properties of the stable process}
\label{sec: stable proc}

\subsection{Stable process, bridge and excursion}

In this part, we gather some well-known facts about stable processes, the associated bridge and excursion processes that will be useful for our proof. Throughout the discussion, we fix the value of a real number $\alpha\in (1, 2)$ and let us recall the spectrally positive $\alpha$-stable process $\rY^{\alpha}=(Y^{\alpha}_{t})_{t\ge 0}$ defined on some probability space $(\Omega, \mathcal F, \mathbb P)$, whose Laplace exponent is given by $\psi(\lambda)=\lambda^{\alpha}$.  
We denote by $p_{t}$ the probability density function of $Y^{\alpha}_{t}$. Note that $p_{t}$ has the Laplace transform:
\[
\int_{\R} e^{-\lambda x}p_{t}(x) dx = \exp(t\lambda^{\alpha}), \quad \lambda\ge 0, t\ge 0.
\] 
In particular, this shows that $x\mapsto p_{t}(x)$ is continuous and $\sup_{x\in \R}p_{t}(x)<\infty$ for fixed $t$. A nice property of $\rY^{\alpha}$ that will be important to us is its {\it invariance by scaling}, namely, for any $c>0$, we have
\begin{equation}
\label{eq: scaling}
\big(c^{-\frac{1}{\alpha}}Y^{\alpha}_{ct}\big)_{ t\ge 0} \eqd \rY.
\end{equation}

\paragraph{Stable bridges.}
The {\it bridge process $\rY^{\br}=(Y^{\br}_{t})_{0\le t \le 1}$ for $\rY^{\alpha}$} has right-continuous sample paths with left-hand limits (i.e.~{\it c\`adl\`ag}), and its law is characterised by two properties: (i) $\mathbb P(Y^{\br}_{1}=0)=1$; (ii) an absolute continuity relationship holds for each $t\in (0, 1)$: if $F$ is a bounded continuous functional defined on the Skorokhod space $\mathbb D([0, t], \R)$, then
\begin{equation}
\label{def: bridge}
\mathbb E\big[F(Y^{\br}_{s}; 0\le s\le t)\big] = \mathbb E\left[F(Y^{\alpha}_{s}; 0\le s\le t)\frac{p_{1-t}(-Y^{\alpha}_{t})}{p_{1}(0)}\right].
\end{equation}
See for instance \cite{Bertoin}. Note that the time reversal property of $\rY^{\alpha}$ and \eqref{def: bridge} together imply the following time reversal property for $\rY^{\br}$: let $\mathrm{\widehat Y}^{\br} =(\widehat Y^{\br}_{t})_{0\le t\le 1}$ be defined by $\widehat Y^{\br}_{t}=-Y^{\br}_{(1-t)-}, 0\le t\le 1$; then we have
\begin{equation}
\label{time-reversal}
\mathrm{\widehat Y}^{\br} \eqd \rY^{\br}. 
\end{equation}

\paragraph{Stable excursions.} Denote by $\mathrm I=(I_{t})_{t\ge 0}$ the infimum process of $\rY$: $I_{t}=\inf_{0\le s\le t}Y^{\alpha}_{s}$. The absence of negative jumps in $\rY^{\alpha}$ means that $\mathrm I$ has continuous sample paths. On the other hand, as the sample paths of $\rY^{\alpha}$ have unbounded variations, $0$ is visited instantaneously by the reflected process $\rY^{\alpha}-\mathrm I$. It follows from the excursion theory of Markov processes (see Chapter IV in \cite{Bertoin}) that $-\mathrm I$ serves as a local time for the excursions of $\rY^{\alpha}-\mathrm I$ away from $0$. Denote by $(g_{i}, d_{i})$, $i\in\N$, the connected components of $\{t>0: Y^{\alpha}_{t}>I_{t}\}$ (in other words, $(g_{i}, d_{i})$'s are the excursion intervals); and define $\mathrm e^{i}$ as the excursion on $(g_{i}, d_{i})$:
\[
e^{i}_{t} = Y^{\alpha}_{t+g_{i}}-I_{t+g_{i}} =Y^{\alpha}_{t+g_{i}}- I_{g_{i}} , \quad 0\le t\le d_{i}-g_{i}.
\]
The excursion theory says that the point measure on $[0, \infty)\times \mathbb D([0, \infty), \R)$: 
\[ 
\sum_{i\in \N}\delta_{(-I_{g_{i}}, \,\mathrm e^{i})}
\]
is a Poisson point process, whose intensity is a $\sigma$-finite measure denoted as $\mathbf N$. The measure $\mathbf N$ is often referred to as the {\it excursion measure}. 
Let us write $\re=(e_{t})_{t\ge 0}$ for the canonical process on the Skorokhod space $\mathbb D([0, \infty), \R)$ and denote by $\zeta=\zeta(\re)=\inf\{ t>0: e_{s}=0\,\forall\,s\ge t\}$ its lifetime. 
As a consequence of the scaling property \eqref{eq: scaling}, 
there exists a probability measure $\nr$ on $\mathbb D([0, \infty), \R)$, called the {\it normalised excursion measure}, so that we can disintegrate $\mathbf N$ with respect to $\zeta$; more precisely, we have for any measurable function $F$ on the Skorokhod space, \begin{equation}
\label{eq: normalise}
\mathbf N\big(F(\re) \big) =\int \mathbf N(\zeta\in dr) \nr\big(F((r^{\frac{1}{\alpha}}e_{t/r})_{t\ge 0} )\big). 
\end{equation}
Intuitively, $\nr$ is the law of an excursion of $\rY^{\alpha}-\mathrm I$ conditioned on $\zeta=1$. Let us also note that the term $\mathbf N(\zeta\in dr)$ in \eqref{eq: normalise} is also known. Indeed, the fluctuation theory applied to $\rY^{\alpha}$ implies that the right-continuous inverse of $\mathrm I$: $T_{x}=\inf\{t>0: I_{t}<-x\}$, $x>0$, is a stable subordinator of index $1/\alpha$. Combined with the excursion theory, this leads to 
\begin{equation}
\label{eq: zeta-dis}
\mathbf N(1-e^{-\lambda \zeta}) = \psi^{-1}(\lambda)= \lambda^{\frac{1}{\alpha}}, \quad \lambda>0.
\end{equation}
Inverting the Laplace transform, we find that $\mathbf N(\zeta\in dr) =(\alpha\Gamma(1-1/\alpha))^{-1} r^{-\frac{1}{\alpha}-1}dr$. 

The Vervaat transformation provides another construction for the normalised excursion. 
Recall the process $\mathrm X^{\alpha}$ from \eqref{def: Vervaat}. Chaumont \cite{Chau97} shows that 
\[
\mathrm X^{\alpha} \text{ under } \mathbb P\  \eqd\  \re \text{ under }\nr\,.
\]  

\subsection{Jumps in a stable excursion}
\label{sec: jumps}

We give a proof of Lemma \ref{lem: normalisation} here, based upon the various properties of stable bridges and excursion processes recalled above. 

\begin{proof}[Proof of Lemma~\ref{lem: normalisation}]
For $x>0$, let us define
\[
\varphi(x)=e^{-x}-1+x \le \min\big\{\tfrac12x^{2}, x\big\}. 
\]
Since the Vervaat transformation \eqref{def: Vervaat} preserves the jump sizes, the first limit in \eqref{eq: asymp} is equivalent to the following 
\begin{equation}
\label{f-limit}
t^{-\alpha}\sum_{0\le u\le 1}\varphi\big(t\Delta Y^{\br}_{u}\big)\xrightarrow{t\to\infty} 1 \quad \text{ in probability}.
\end{equation}
The first step in confirming \eqref{f-limit} consists in showing that for all $s>0$, we have 
\begin{equation}
\label{f-limit'}
t^{-\alpha}\sum_{0\le u\le s}\varphi\big(t\Delta Y^{\alpha}_{u}\big) \xrightarrow{t\to\infty} s \quad \text{ in probability}.
\end{equation}
Note that we only need to consider jumps of magnitudes $\le 1$, since there is only a finite number of jumps greater than $1$ on $[0, s]$ and $\varphi(x)\le x$.  On the other hand, $(\Delta Y^{\alpha}_{u})_{u\ge 0}$ has the distribution of a Poisson point process of intensity $\pi(dx)=c_{\alpha}x^{-1-\alpha}\mathbf 1_{\{x>0\}}dx$, with $c_{\alpha}= \alpha(\alpha-1)/\Gamma(2-\alpha)$. The compensation formula for the Poisson point process yields that 
\[
\mathbb E\Big[\sum_{0\le u\le s}\varphi\big(t\Delta Y^{\alpha}_{u}\big)\mathbf 1_{\{\Delta Y^{\alpha}_{u}\le 1\}}\Big]=s\int_{(0, 1)} \varphi(tx)\pi(dx) =st^{\alpha}-\int_{(1, \infty)}\varphi(tx) \pi(dx).
\]
Using the bound $\varphi(tx)\le tx$, we deduce that the second term above is at most $\frac{c_{\alpha}}{\alpha-1}st$. It follows that
\begin{equation}
\label{cal: epec}
t^{-\alpha}\,\mathbb E\Big[\sum_{0\le u\le s}\varphi\big(t\Delta Y^{\alpha}_{u}\big)\mathbf 1_{\{\Delta Y^{\alpha}_{u}\le 1\}}\Big]\to s, \ \text{ as } t\to\infty.
\end{equation}
We also deduce from the exponential formula for the Poisson point process that
\[
\Var\Big[\sum_{0\le u\le s}\varphi\big(t\Delta Y^{\alpha}_{u}\big)\mathbf 1_{\{\Delta Y^{\alpha}_{u}\le 1\}}\Big] = s\int_{(0, 1)} \varphi(tx)^{2}\pi(dx) \le st^{2}\int_{(0, 1)}x^{2}\pi(dx)
\le \frac{c_{\alpha}}{2-\alpha} st^{2}.
\]
%
%
%
%
Together with \eqref{cal: epec} and Markov's inequality, this implies \eqref{f-limit'} for each fixed $s>0$. 
Now take $s\in (1-\delta, 1)$. On the one hand, \eqref{def: bridge} together with \eqref{f-limit'} implies that
\begin{equation}
\label{int: s'}
\mathbb P\Big(\Big| t^{-\alpha}\sum_{0\le u\le s}\varphi\big(t\Delta Y^{\br}_{u}\big)-s\Big| >\delta \Big) \le \mathbb P\Big(\Big|t^{-\alpha}\sum_{0\le u\le s}\varphi\big(t\Delta Y^{\alpha}_{u}\big)-s\Big| >\delta\Big)\cdot\frac{\sup_{x}p_{1-s}(x)}{p_{1}(0)}\overset{t\to\infty}{\longrightarrow} 0. 
\end{equation}
On the other hand,  it follows from the time reversal property \eqref{time-reversal} that 
\begin{align*}
\mathbb P\Big(\sum_{s\le u\le 1}\varphi\big(t\Delta Y^{\br}_{u}\big)>\sqrt\delta  t^{\alpha}\Big) &= \mathbb P\Big(\sum_{0\le u\le 1-s}\varphi\big(t\Delta Y^{\br}_{u}\big)>\sqrt\delta t^{\alpha}\Big)\\
& \le \mathbb P\Big(\sum_{0\le u\le 1-s}\varphi\big(t\Delta Y^{\alpha}_{u}\big)>\sqrt\delta t^{\alpha}\Big) \cdot\frac{\sup_{x}p_{s}(x)}{p_{1}(0)} \\
& \le \frac{\mathbb E[\sum_{0\le u\le 1-s}\varphi\big(t\Delta Y^{\alpha}_{u}\big)]}{\sqrt\delta t^{\alpha}} \cdot \frac{\sup_{x}p_{s}(x)}{p_{1}(0)}\le \sqrt\delta \cdot \frac{\sup_{x}p_{s}(x)}{p_{1}(0)},
\end{align*}
where we have used \eqref{def: bridge} in the first inequality, Markov's inequality in the second, and the compensation formula in the third. 
Combining the above with \eqref{int: s'}, we deduce the convergence in \eqref{f-limit} 
by first taking $t\to\infty$ and then $\delta\to 0$. The second limit in \eqref{eq: asymp} readily follows from the first, as $\Psi^{-1}_{\bDelta^{\downarrow}_{\alpha}}$ is the inverse function of $\Psi_{\bDelta^{\downarrow}_{\alpha}}$. For the third one, we note that 
once again the Vervaat transformation and the finite numbers of large jumps reduce the  proof to the following:
\[
\epsilon^{\alpha-1}\sum_{0\le u\le 1}\Delta Y^{\br}_{u}\mathbf 1_{\{\epsilon<\Delta Y^{\br}_{u}\le 1\}} \xrightarrow{\epsilon\to 0} \frac{\alpha}{\Gamma(2-\alpha)} \quad \text{ in probability}.
\]
Its proof is quite similar to that of \eqref{f-limit}: it suffices to replace \eqref{cal: epec} with 
\[
\epsilon^{\alpha-1}\mathbb E\Big[\sum_{0\le u\le s}\Delta Y^{\alpha}_{u} \mathbf 1_{\{\epsilon<\Delta Y^{\alpha}_{u}\le 1\}}\Big]=\frac{\alpha}{\Gamma(2-\alpha)}\cdot s(1-\epsilon^{\alpha-1}) \xrightarrow{\epsilon\to 0} \frac{\alpha s}{\Gamma(2-\alpha)}, 
\] 
and the bound on the variance with 
\[
\Var\Big[\sum_{0\le u\le s}\Delta Y^{\alpha}_{u} \mathbf 1_{\{\epsilon<\Delta Y^{\alpha}_{u}\le 1\}}\Big] = \frac{\alpha(\alpha-1)}{\Gamma(3-\alpha)}\cdot s (1-\epsilon^{2-\alpha}).
\]
We therefore omit the detail. 
\end{proof}

\section{Real trees and stable trees}

\subsection{Real trees and distances between metric spaces}
\label{sec: GH}

This subsection is a recap on real trees, their encodings by real-valued functions and the Gromov--Hausdorff topology. 

A {\it real tree} $(T, d)$ is a complete metric space which satisfies the following two properties for all pairs $(x, y)$ of points of $T$. First, there is a geodesic connecting $x$ to $y$, namely there is an isometric embedding $f: [0, d(x, y)]\to T$ so that $f(0)=x$ and $f(d(x, y))=y$; in the sequel, we will denote by $\llb x, y\rrb_{T}=f([0, d(x, y)])$ this geodesic. 
Second, the aforementioned geodesic provides the unique path between $x$ and $y$; more precisely, 
if $g: [0, 1]\to T$ is a continuous mapping with $g(0)=x$ and $g(1)=y$, then necessarily $g([0, 1])=\llb x,y\rrb_{T}$. 

We note that the above definition of real tree is an extension to our concept of a (graph) tree as a connected and loop-free graph, where the length of the unique path between two vertices determines their graph distance. In particular, if we take a finite graph tree and replace each of its edges by the $[0, 1]$ interval, this will give us a somewhat boring example of real trees. More exciting examples can be obtained with the help of stochastic processes. To that end, let us first recall how to extract real trees from continuous excursion-like functions. 

Throughout this subsection, let $f$ be a {\it continuous} real-valued function with compact support. Denote by $\zeta=\zeta(f)=\inf\{t>0: f(s)=0 \,\forall\, s\ge t\}$ to be the upper end of its support, or simply its {\it lifetime}. We further suppose that $f(0)=f(\zeta)=0$ and $f(t)>0$ for all $t\in (0, \zeta)$. We introduce the following symmetric function on $[0, \zeta]$:
\begin{equation}
\label{def: f-encoding}
d_{f}(s, t) = f(s) + f(t) -2m(f, s,t), \quad \text{where} \quad m(f, s, t) = \inf\{f(u): s\wedge t \le u\le s\vee t\}.
\end{equation}
It turns out that $d_{f}$ verifies the triangle inequality. To turn it into a genuine metric, we introduce the equivalence relation $\sim_{f}$ on $[0, \zeta]$: we say $s\sim_{f}t $ if and only if $d_{f}(s, t)=0$. 
Let $T_{f}=[0, 1]/\sim_{f}$ be the quotient space; then $d_{f}$ defines a metric on it. Moreover, the pair $(T_{f}, d_{f})$ is a compact real tree \cite{DuLeG05}. 

Comparing \eqref{def: f-encoding} with the definition \eqref{def: H-encoding} of the stable tree, we see that the $\alpha$-stable tree $(\cT_{\alpha}, d_{\alpha})$ is the real tree ``extracted'' from the height process $\mathrm H$. We wish to consider $(\cT_{\alpha}, d_{\alpha})$ as a ``random real tree''; this is possible as we will shortly see that the space of compact real trees is a Polish space under the so-called Gromov--Hausdorff topology. 

If $(X, d_{X}), (Y, d_{Y})$ are two {\it compact} metric space, their mutual {\it Gromov--Hausdorff distance} is defined as
\[
\dgh\big((X, d_{X}), (Y, d_{Y})\big) = \inf \dhau\big(\phi(X), \varphi(Y)\big), 
\]
where the infimum is over all the isometric embeddings $\phi: X\to Z$ and $\varphi: Y\to Z$ into a common metric space $(Z, d_{Z})$, and $\dhau$ is the Hausdorff distance on the compacts sets of $Z$. In particular, two compact metric spaces are {\it isometric} if their Gromov--Hausdorff distance is null. Denote by $\mathbb T_{c}$ the set of isometry equivalence classes of compact real trees. Then $(\mathbb T_{c}, \dgh)$ is a Polish space (\cite{Evans}). 

Real trees such as stable trees considered in this paper are rooted and equipped with a probability measure. We can refine the notion of Gromov--Hausdorff distance to take into account these additional features. Let $(X, d_{X}), (Y, d_{Y})$ be as before. Suppose that $x\in X, y\in Y$, and $\mu_{X}, \mu_{Y}$ are respectively (Borel) probability measures on $X$ and $Y$. Then the {\it pointed Gromov--Hausdorff--Prokhorov distance} between $X$ and $Y$ is given by
\begin{multline*}
\dghp\big((X, d_{X}, x, \mu_{X}), (Y, d_{Y}, y, \mu_{Y})\big) \\
= \inf \big\{\dhau\big(\phi(X), \varphi(Y)\big)+d_{Z}(\phi(x), \varphi(y))+\dpr(\mu_{X}\circ\phi^{-1}, \mu_{Y}\circ\varphi^{-1})\big\},
\end{multline*}
where as before the infimum is over all the isometric embeddings $\phi: X\to Z$ and $\varphi: Y\to Z$ into a common metric space $(Z, d_{Z})$, and $\dpr$ is the Prokhorov distance for probability measures on $Z$. 

Equipping a real tree with a probability measure on a real tree not only facilitates its analysis, but also has measure theoretic implications, as explained in what follows. The inhomogeneous continuum random trees are not all compact: some are merely complete as metric space (see \cite{AP02}).  
We will call a complete metric spaces equipped with a Borel probability measure as a {\it measured metric space}. 
For two measured metric spaces $(X, d_{X}, \mu_{X})$ and $(Y, d_{Y}, \mu_{Y})$, their {\it Gromov--Prokhorov} distance is defined as
\[
\dgp\big((X, d_{X}, \mu_{X}), (Y, d_{Y}, \mu_{Y})\big) 
= \inf \dpr(\mu_{X}\circ\phi^{-1}, \mu_{Y}\circ\varphi^{-1}),
\]
where the infimum is over all the isometric embeddings $\phi: \supp(\mu_{X})\to Z$ and $\varphi: \supp(\mu_{Y})\to Z$ into a common metric space $(Z, d_{Z})$, with $\supp(\mu_{X}), \supp(\mu_{Y})$ standing for the respective support sets of $\mu_{X}, \mu_{Y}$. 
Two measured metric spaces are {\it equivalent} if their Gromov--Prokhorov distance is null. Denote by $\mathbb T$ the set of all equivalence classes of measured metric spaces that are also real trees. Then $\mathbb T$ is a Polish space under the topology induced by $\dgp$ (\cite{GPW09}).

\subsection{Recovery of the stable tree}
\label{sec: stable-recov}

This section contains the proof of Proposition~\ref{prop: stable-recov}. 
Recall the i.i.d.~sequence of points $(\sigma_{i})_{i\ge 1}$; each $\sigma_{i}$ is a leaf of $\cT_{\alpha}$ as $\mu_{\alpha}$ only charges $\Lf(\cT_{\alpha})$. Recall also that 
\[
\cT_{k}=\bigcup_{1\le i\le k}\llb r_{\alpha}, \sigma_{i}\rrb_{\cT_{\alpha}}
\]
 is the subtree spanned by $\sigma_{1}, \sigma_{2}, \cdots, \sigma_{k}$ and the root $r_{\alpha}$.  
Let $\mathcal V_{k}=\Br(\cT_{k})\cup\Lf(\cT_{k})\cup\{r_{\alpha}\}$, which is a finite set. 
We define $T_{k}=(V(T_{k}), E(T_{k}))$, $k\ge 1$ to be a sequence of (graph) trees that satisfy the following properties:
\begin{itemize}
\item
{\bf $T_{k}$ has the shape of $\cT_{k}$}: there is a bijection $f_{k}:  \mathcal V_{k}\to V(T_{k})$ such that 
\begin{equation}
\label{def: shape}
\{f_{k}(x), f_{k}(y)\}\in E(T_{k}) \quad\text{ if and only if }  \quad\rrb x, y\llb_{\cT_{\alpha}} \,\cap\, \mathcal V_{k}=\varnothing;
\end{equation}
\item
{\bf the labelling is consistent across $k$}: $f_{k}$ is a restriction of $f_{k+1}$ to $\mathcal V_{k}$, $k\in \N$.  
\end{itemize}
%
%
It is not difficult to see that up to a choice in vertex labelling, the sequence $(T_{k})_{k\in \N}$ exists in a unique way. Note also from \eqref{def: shape} that we have
\begin{equation}
\label{id: deg}
\deg(x, \cT_{k})=\deg(f_{k}(x), T_{k}), \quad \forall\, x\in \mathcal V_{k}, k\ge 1.
\end{equation}
Therefore, the statements in Proposition~\ref{prop: stable-recov} will follow from the following properties of $\cT_{k}$ and $\cT_{\alpha}$. We recall the convention that $\deg(v, \cT_{k})=-\infty$ if $v\notin \cT_{k}$. 
\begin{prop}
\label{prop: stable-recov'}
The following statements  hold true $\mathbb P$-a.s.
\begin{enumerate}[(i)]
\item
For each $v\in \Br(\cT_{\alpha})$\,, we have
\begin{equation}
\label{def: loc-stable''}
\Delta^{\alpha}(v)=\lim_{k\to\infty} \frac{\deg(v, \cT_{k})}{k^{1/\alpha}},
\end{equation}
where the limit exists almost surely. 
\item
For $v, v'\in \cT_{k}$\,, $k\ge 1$, we have 
\begin{equation}
\label{def: dist-stable'}
d_{\alpha}(v, v') =\frac{\Gamma(2-\alpha)}{\alpha}\cdot \lim_{\epsilon\to 0+} \epsilon^{\alpha-1}\#\big\{w\in \Br(\cT_{\alpha})\cap\llb v, v'\rrb_{\cT_{\alpha}}: \Delta^{\alpha}(w)>\epsilon \big\},
\end{equation}
where the above limit exists in probability. Moreover, $\Br(\cT_{\alpha})$ is dense everywhere in $(\cT_{\alpha}, d_{\alpha})$. 
\item
The sequence of probability measures $\frac1k\sum_{1\le i\le k}\delta_{\sigma_{i}}$, $k\ge 1$, converges a.s.~to $\mu_{\alpha}$ in the weak topology of $(\cT_{\alpha}, d_{\alpha})$. 
\end{enumerate}
\end{prop}

\begin{proof}[Proof of Proposition \ref{prop: stable-recov'}]
The statement in {\it (i)} is undoubtedly a well-accepted fact about stable trees; however we have failed to find a reference. So we provide a proof of it in Appendix \ref{sec: A1}, relying upon the Poissonian marking technique used in \cite{DuLG02}. 
To prove {\it (ii)}, let us first suppose $v'$ to be the root $r_{\alpha}$. 
It follows from Lemma 5.1 in \cite{BrDuWa22+} that  for $t\in (0,1)$ and  $s\in [0, t)$ with $\Delta X^{\alpha}_{s}>0$, we have $p(s)\in \llb r_{\alpha}, p(t)\rrb $ if and only if $X^{\alpha}_{s-}< \inf_{s\le u\le t}X^{\alpha}_{u}$. Together with \eqref{eq: corres},  this proves \eqref{def: dist-stable'} in the case that $v'=r_{\alpha}$. The general case readily follows since 
\begin{equation}
\label{id: tree}
d_{\alpha}(v, v')=d_{\alpha}(r_{\alpha}, v)+d_{\alpha}(r_{\alpha}, v')-2\cdot d_{\alpha}(r_{\alpha}, v\wedge v'),
\end{equation}
where $v\wedge v'$ is the most recent common ancestor of $v$ and $v'$. 
Next, since  $(\cT_{\alpha}, d_{\alpha})$ is the image of $[0, 1]$ by the continuous mapping $p$, which maps the jump times of $\bX^{\alpha}$ to $\Br(\cT_{\alpha})$,  
the fact that $\Br(\cT_{\alpha})$ is everywhere dense in $\cT_{\alpha}$ readily follows from the fact that the jump times of $\bX^{\alpha}$ are everywhere dense in $[0, 1]$. Finally, {\it (iii)} follows from the Glivenko--Cantelli Theorem after conditioning on $(\cT_{\alpha}, d_{\alpha})$. 
\end{proof}

%

\begin{proof}[Proof of Proposition~\ref{prop: stable-recov}]
For {\it (i)}, if $v$ is a branch point of some $T_{k}$, then according to \eqref{id: deg}, $f_{k}^{-1}(v)\in \Br(\cT_{k})\subset \Br(\cT_{\alpha})$. In that case, the limit in \eqref{def: loc-stable'} follows from  \eqref{def: loc-stable''} and \eqref{id: deg}, and we have $\tilde\Delta^{\alpha}(v)=\Delta^{\alpha}(f_{k}^{-1}(v))$. If, on the other hand, $v$ is never a branch point in the sequence $(T_{k})_{k\ge 1}$, then $\deg(v, T_{k}), k\ge 1$, is bounded and therefore $\widetilde{\Delta}^{\alpha}(v)=0$ almost surely. Together with \eqref{eq: corres}, this completes the proof of {\it (i)}. For {\it (ii)}, thanks to {\it (i)} and \eqref{def: shape}, 
for $k\ge 1$ and $v, v'\in \mathcal V_{k}$, we have $N^{\alpha}_{\epsilon}(f_{k}(v), f_{k}(v'))=\#\{ w\in \Br(\cT_{\alpha})\cap\llb v, v'\rrb_{\cT_{k}}: \Delta^{\alpha}(w)>\epsilon\}$, so that the existence of the limit in \eqref{def: dist-stable} follows from \eqref{def: dist-stable'}, and we have $\widetilde d_{\alpha}(f_{k}(v), f_{k}(v'))=d_{\alpha}(v, v')$ a.s. Since $(\cT_{\alpha}, d_{\alpha})$ is the completion of $(\cup_{k}\mathcal V_{k}, d_{\alpha})$, the rest of the statements in {\it (ii)} follow. Finally, as $\nu^{\alpha}_{k}$ is the image of $\frac1k\sum_{1\le i\le k}\delta_{\sigma_{i}}$ by $f_{k}$, {\it (iii)} also holds true. 
\end{proof}

\section{Recovery of the ICRT}
\label{sec: icrt-recov}

Let us recall that $\cR'_{k}$ is the subtree of $\cT$ spanning the i.i.d.~points $\eta'_{1}, \eta'_{2}, \dots, \eta'_{k}$. 
From the Line-breaking construction, it is not difficult to see that all the branch points of $\cT$ are given by the images of $\xi_{i, 1}, i\ge 1$, which we still denote as $\xi_{i, 1}$.%
%
We can take the same steps as in the stable case to define an increasing sequence of discrete trees $(R'_{k})_{k\ge 1}$ which represent the shapes of $(\cR'_{k})_{k\ge 1}$. Proposition~\ref{prop: icrt-recov} will then be a consequence of the corresponding properties of $(\cR'_{k})_{k\ge 1}$, which have been mostly proved in \cite{BHW22+}. 

\begin{prop}
\label{prop: icrt-recov'}
Let the functions $\Psi_{\bth}$, $\Psi^{-1}_{\bth}$ and $\gamma_{\bth}$ be defined as in Proposition~\ref{prop: icrt-recov}. 
The following statements  hold true $\mathbf P^{\bth}$-a.s.
\begin{enumerate}[(i)]
\item
For each $i\in \N$, we have
\begin{equation}
\label{def: loc-icrt''}
\theta_{i}=\lim_{k\to\infty} \frac{\deg(\xi_{i, 1}, \cR'_{k})}{\Psi^{-1}_{\bth}(k)},
\end{equation}
where the limit exists in probability. 
\item
For  $v, v'\in \cR'_{k}$\,, $k\ge 1$, we have 
\begin{equation}
\label{def: dist-icrt'}
d(v, v') = \lim_{\epsilon\to 0+} \frac{1}{\gamma_{\bth}(\epsilon)}\#\big\{i\in \N: \theta_{i}>\epsilon \text{ and } \xi_{i, 1}\in \llb v, v'\rrb_{\cT}\big\},
\end{equation}
where the above limit exists in probability. Moreover, $\Br(\cT)$ is dense everywhere in $(\cT, d)$.
\end{enumerate}
\end{prop}

\begin{proof}[Proof of Proposition~\ref{prop: icrt-recov'}]
The statements in {\it (i)} correspond to Proposition 1 in \cite{BHW22+}. For the limit in \eqref{def: dist-icrt'}, the arguments are based upon the proof of Proposition 2 in \cite{BHW22+}. 
Let $\ell\in \N$ and denote by $\mathrm m^{(\ell)}=\sum_{1\le i\le \ell}\theta_{i}\delta_{\xi_{i, 1}}$, a finite measure on $\cT$. Proposition 5(b) of \cite{Al00} implies that for each $k, \ell\in \N$, 
\[
\Big(\cR'_{k}, \mathrm m^{(\ell)}(\cdot \cap \cR'_{k})\Big) \text{ has the same distribution as } \Big(\cR_{k}, \mathrm m^{(\ell)}(\cdot \cap \cR_{k})\Big)\ \text{ under } \mathbb P^{\bth}.
\]
So it suffices to prove \eqref{def: dist-icrt'} for $v, v'\in \cR_{k}$. 
Thanks to an analogue of \eqref{id: tree} in the ICRT case, we further reduces the case under consideration to $v'=r$ and $v\in \llb r, \eta_{k}\rrb$, $k\ge 1$. 
But the law of $(\eta_{i})_{i\ge 1}$ is exchangeable. Therefore, we only need to consider the case $k=1$. From the Line-breaking algorithm, the branch $\llb r, \eta_{1}\rrb$ is simply the image of $[0, \eta_{1}]$ in $\cT$. Therefore, \eqref{def: dist-icrt'} will be a consequence of the following statement:
\begin{equation}
\label{eq: dist-sup}
\sup_{0\le x\le \eta_{1}} \left|\frac{1}{\gamma_{\bth}(\epsilon)}\#\Big\{i\in \N: \theta_{i}>\epsilon \text{ and } \xi_{i, 1}\le x\Big\}-x\right| \overset{\epsilon\to 0}{\longrightarrow} 0, \quad \text{in probability}.
\end{equation}
Let us show \eqref{eq: dist-sup}. It is clear from the Line-breaking algorithm that $(\xi_{i, 1})_{i\ge 1}$ is a collection of independent exponential variables with $\mathbb E[\xi_{i, 1}]=\theta_{i}^{-1}$ and  $\eta_{1}=\min\{\xi_{i, 2}: i\ge 1\}$. For $t\ge 0$ and $\epsilon>0$, we define
\[
L_{\epsilon}(t)= \sum_{i: \theta_{i}>\epsilon} \mathbf 1_{\{\xi_{i, 1}\le t\}} \quad \text{and}\quad  M_{i}(t)=\mathbf 1_{\{\xi_{i, 1}\le t\}}-\theta_{i}(t\wedge \xi_{i, 1}), \  \text{ for each } i\ge 1. 
\]
Note that $(M_{i}(t))_{t\ge 0}$ is a martingale with respect to the natural filtration of $(\mathbf 1_{\{\xi_{i, 1}\le t\}})_{t\ge 0}$ and that $\mathbb E[M_{i}(t)^{2}]\le \mathbb E[\mathbf 1_{\{\xi_{i, 1}\le t\}}]=1-\exp(-\theta_{i}t)$. 
Now let
\[
A_{\epsilon}(t)=\sum_{i: \theta_{i}>\epsilon}\theta_{i}(t \wedge \xi_{i, 1}) \quad \text{and}\quad \mathcal M_{\epsilon}(t)=\sum_{i: \theta_{i}>\epsilon}M_{i}(t), \quad t\ge 0. 
\]
Then $(\mathcal M_{\epsilon}(t))_{t\ge 0}$ is a martingale with respect to the natural filtration of $\{(\mathbf 1_{\{\xi_{i, 1}\le t\}})_{t\ge 0}: i\ge 1\}$. 
Thanks to Doob's maximal inequality and the fact that $M_{i}(t)$, $i\ge 1$, are independent, we deduce that
\begin{equation}
\label{bd: mart}
\frac{1}{\gamma_{\bth}(\epsilon)^{2}}\mathbb E\Big[\sup_{s\le t}\mathcal M_{\epsilon}(s)^{2}\Big]\le \frac{4}{\gamma_{\bth}(\epsilon)^{2}}\mathbb E\Big[\mathcal M_{\epsilon}(t)^{2}\Big]\le \frac{4}{\gamma_{\bth}(\epsilon)^{2}}\sum_{i: \theta_{i}>\epsilon}(1-e^{-\theta_{i}t}) \le \frac{4t}{\gamma_{\bth}(\epsilon)}.
\end{equation}
On the other hand, we have
\[
\bigg|\frac{1}{\gamma_{\bth}(\epsilon)}A_{\epsilon}(t)-t\bigg|
=\bigg|\frac{1}{\gamma_{\bth}(\epsilon)}\sum_{i: \theta_{i}>\epsilon}\theta_{i}\int_{0}^{t}\big(\mathbf 1_{\{\xi_{i, 1}> s\}}-1\big)ds\bigg|
=\frac{1}{\gamma_{\bth}(\epsilon)}\sum_{i: \theta_{i}>\epsilon}\int_{0}^{t}\theta_{i}\mathbf 1_{\{\xi_{i, 1}\le s\}}ds,
\]
which is clearly increasing in $t$. It follows that
\begin{equation}
\label{bd: comp}
\mathbb E\bigg[\sup_{s\le t}\bigg|\frac{1}{\gamma_{\bth}(\epsilon)}A_{\epsilon}(s)-s\bigg|\bigg] = \mathbb E\bigg[\frac{1}{\gamma_{\bth}(\epsilon)}\sum_{i: \theta_{i}>\epsilon}\int_{0}^{t}\theta_{i}\mathbf 1_{\{\xi_{i, 1}\le s\}}ds\bigg] \le \frac{\|\bth\|^{2}t^{2}}{\gamma_{\bth}(\epsilon)}. 
\end{equation}
Note that $L_{\epsilon}(t)=A_{\epsilon}(t)+\mathcal M_{\epsilon}(t)$. Therefore, \eqref{bd: mart} and \eqref{bd: comp} yield that
\begin{equation}\label{cv: L}
\mathbb E\bigg[\sup_{s\le t}\bigg|\frac{L_{\epsilon}(s)}{\gamma_{\bth}(\epsilon)}-s\bigg|\bigg] \to 0, \quad \forall\, t>0. 
\end{equation}
As the law of $\eta_{1}$ is tight: $\mathbb P(\eta_{1}>t)\to 0$ as $t\to\infty$. We deduce from this and \eqref{cv: L} the convergence in \eqref{eq: dist-sup}. Since $\Br(\cT)=\{\xi_{i, 1}: i\in \N\}$, it is dense in $\cT$ as $\sum_{i}\theta_{i}=\infty$. 
This completes the proof. 
\end{proof}

\begin{proof}[Proof of Proposition~\ref{prop: icrt-recov}]
There is an obvious correspondence between the branch points of $R_{k}$ and $\Br(\cR_{k})$. On the other hand, we have seen $\Br(\cT)=\cup_{k\ge 1}\Br(\cR_{k})=\{\xi_{i, 1}: i\in \N\}$. Therefore, a branch point of $R_{k}$ must correspond to some 
$\xi_{i, 1}$. Then point {\it(i)} of Proposition~\ref{prop: icrt-recov'} says that in that case its local time is given by $\theta_{i}$. 
Rest of the arguments are similar to the ones found in the proof of Proposition~\ref{prop: stable-recov}, and are therefore omitted.\end{proof}

\section{Trees embedded in  c\`adl\`ag functions}
\label{sec: coding}

Throughout this section, we suppose that $\mathbf x=(x(t))_{t\ge 0}\in \mathbb D(\R_{+}, \R)$ satisfying:
\begin{itemize}
\item
Finite support: $\zeta(\mathbf x)=\inf\{t: x(s)=0 \ \forall\, s\ge t\}\in [0, \infty)$;
\item
Positive values: $x(t)\ge 0$ for all $0\le t\le \zeta(\mathbf x)$.
\item
Positive jumps: $\Delta x(t)=x(t)-x(t-)\ge 0$ for all $0\le t\le \zeta(\mathbf x)$. 
\end{itemize}
We will also need the following notation. For $s\in [0, \zeta(\mathbf x))$, denote 
\begin{equation}
\label{def: sigma-time}
\sigma_{\mathbf x}(s)=\inf\big\{t>s: x(t)<x(s-)\big\}\in [s, +\infty].
\end{equation}
We observe that for $0\le s_{1}<s_{2}<\zeta(\mathbf x)$, we have 
\[
\text{either }\  \big(s_{1}, \sigma_{\mathbf x}(s_{1})\big)\cap \big(s_{2}, \sigma_{\mathbf x}(s_{2})\big)=\varnothing \ \text{ or } \ \big(s_{2}, \sigma_{\mathbf x}(s_{2})\big)\subseteq \big(s_{1}, \sigma_{\mathbf x}(s_{1})\big).
\] 
Our aim here is to formalise a notion of genealogy on the set $[0, \zeta(\mathbf x)]$, in which points of $(s, \sigma_{\mathbf x}(s))$ are  descendants of $s$. More precisely, let $\mathbf u_{k}=\{u_{1}, u_{2}, \dots, u_{k}\}$ be a collection of $k$ distinct points of $[0, \zeta(\mathbf x)]$. We will define a (discrete) tree $\mathscr T(\mathbf x; \mathbf u_{k})$ as a function of $\mathbf x$ and $\mathbf u_{k}$.  In the case that $\mathbf x$ has bounded variations, we will see that the genealogy coincides with the one induced by the LIFO construction in Le Gall--Le Jan \cite{LGLJ98}. When $\mathbf x$ is either $\mathrm X^{\alpha}$ or $\mathrm X^{\bth}$ and the $u_{i}$'s are uniformly distributed,  we will show that $\mathscr T(\mathbf x; \mathbf u_{k})$ has the distribution of the $k$-leafed spanning trees of respectively the stable tree and the $\bth$-ICRT. 

\smallskip
\noindent
{\bf Ordered rooted tree. }
For the definition of  $\mathscr T(\mathbf x; \mathbf u_{k})$, it will be convenient to work with ordered rooted trees. So let us first recall Neveu's formalism for these trees. Let $\mathbb U=\{\varnothing\}\cup\bigcup_{n\ge 1}\mathbb N^{n}$. A finite subset $t\subset \mathbb U$ is an {\it ordered rooted tree} if it satisfies: (a) $\varnothing\in t$; (b) if $v=(v_{1}, v_{2}, \dots, v_{n-1}, v_{n})\in t$, then $(v_{1}, v_{2}, \cdots, v_{n-1})\in t$; we call $(v_{1}, v_{2}, \cdots, v_{n-1})$ the {\it parent} of $v$;  (c) for all $v=(v_{1}, v_{2}, \dots, v_{n})\in t$, there is some integer $k\ge 0$ so that $(v_{1}, v_{2}, \cdots, v_{n}, i)\in t$ if and only if $i\le k$. 
An ordered rooted tree can be built by taking a finite sequence of ordered rooted trees and then gluing them to a common root. Formally, for $v=(v_{1}, v_{2}, \dots, v_{n})\in \mathbb U$, we introduce the shift operator $\theta_{v}: \mathbb U\to \mathbb U$ as $\theta_{v}(w)=(v_{1}, v_{2}, \dots, v_{n}, w_{1}, w_{2}, \cdots, w_{k})$ if $w=(w_{1}, w_{2}, \dots, w_{k})$. 
If $p\in \mathbb N$ and $t_{1}, t_{2}, \dots, t_{p}$ are ordered rooted trees, then 
\[
t:=\{\varnothing\}\cup\bigcup_{i=1}^{p} \theta_{i}(t_{i}) \quad \text{with} \quad \theta_{i}(t_{i}):=\{\theta_{(i)}(u): u\in t_{i}\}
\]
is also an ordered rooted tree. Graphically speaking, the tree $t$ is formed by connecting the roots of $t_{i}$ to a common root $\varnothing$. 
We will refer to $t_{i}, 1\le i\le p$, as the {\it subtrees above the root} in $t$. 

\begin{figure}
\centering
\includegraphics[height = 12cm]{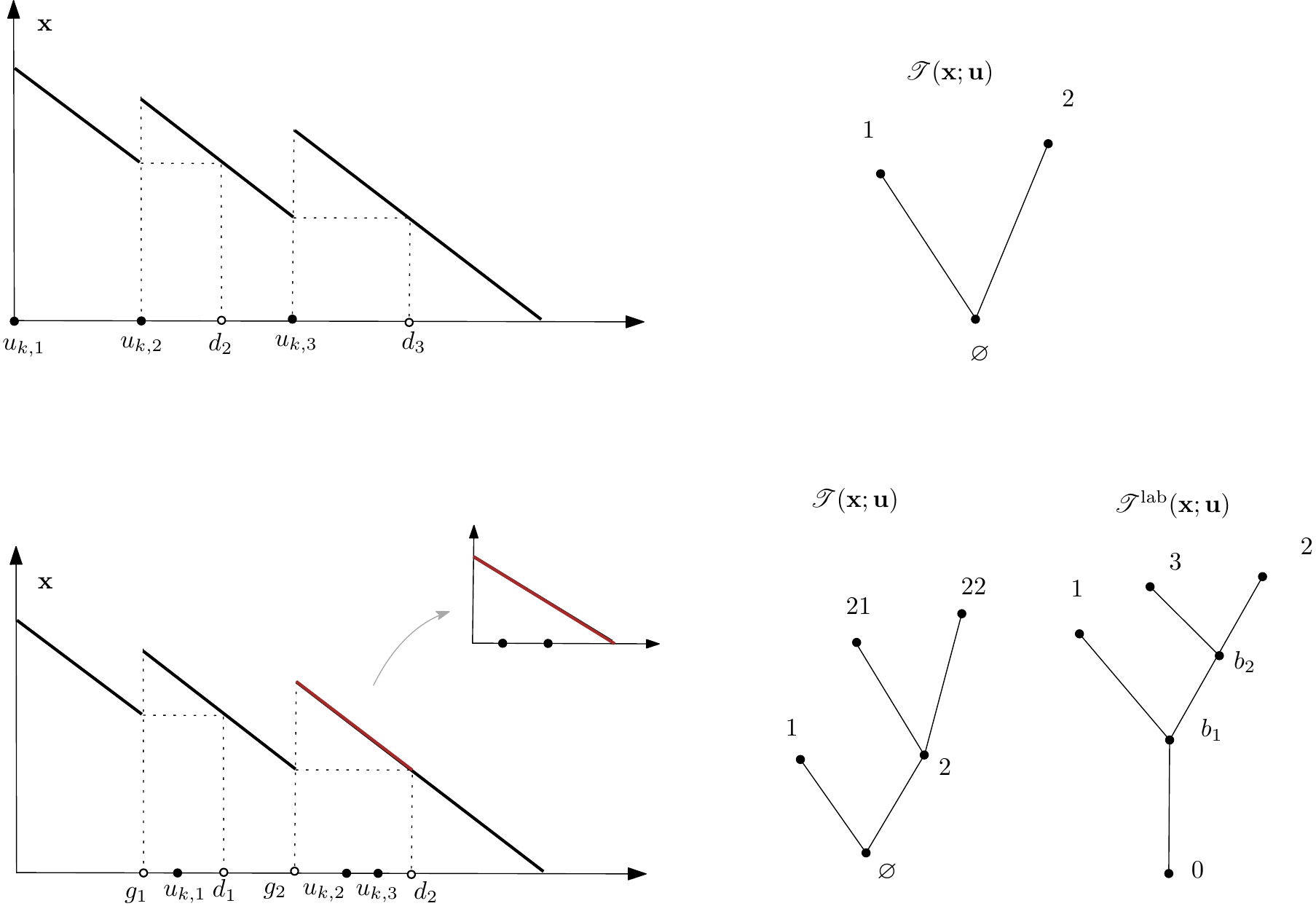}
\caption{\label{fig: spanning}Two examples of $\mathscr T(\mathbf x; \mathbf u)$. In the upper line, $u_{i}$'s coincide with the three jump times. In this example, $b=0$, $\sigma_{\mathbf x}(b)=\zeta(x)$, and $\mathbf z=\mathbf x$. We also have $g_{1}=0, d_{1}=\sigma_{\mathbf x}$, $g_{2}=u_{k, 2}$ and $g_{3}=u_{k, 3}$. Right to it, a depiction of the corresponding $\mathscr T(\mathbf x; \mathbf u)$. 
In the lower line, after the first generation is found, the construction is applied to the red segment to build the second generation of the tree. 
The corresponding  $\mathscr T(\mathbf x; \mathbf u)$ and the labelled version $\mathscr T^{\mathrm{lab}}(\mathbf x; \mathbf u)$ are given on the right. }
\end{figure}

\smallskip
\noindent
{\bf Definition of $\mathscr T(\mathbf x; \mathbf u_{k})$. }
If $k=1$ or $\zeta(\mathbf x)=0$, then $\mathscr T(\mathbf x; \mathbf u_{k})=\{\varnothing\}$. 
If $k\ge 2$ and $\zeta(\mathbf x)>0$, let us denote by $u_{k, 1}< u_{k, 2}< \cdots < u_{k, k}$ the re-arrangement of $\mathbf u_{k}$ in increasing order. 
We will need the following notation: for $0\le s\le t\le \zeta(\mathbf x)$ and $r\ge 0$, let 
\[
m(\mathbf x, s, t)=\inf_{u\in [s, t]}x(u) \quad\text{and}\quad \tau(\mathbf x, t, r) =\inf\{s\le t: m(\mathbf x, s, t) \ge r\}, 
\]
with the convention that $\inf\varnothing = \infty$. 
Let us set 
\begin{equation}
\label{def: mrca}
b=\tau\big(\mathbf x, u_{k, 1}, m(\mathbf x, u_{k, 1}, u_{k, k})\big)=\inf\Big\{s\le u_{k, 1}: \inf_{u\in [s, u_{k, 1}]}x(u)\ge \inf_{u\in [u_{k, 1}, u_{k, k}]}x(u)\Big\}, 
\end{equation}
which will serve as the most recent common ancestor of the $u_{i}$'s. Note that we always have $b<\infty$ and as a matter of fact $b\le u_{k, 1}$. Moreover, $x(b-)\le m(\mathbf x, u_{k, 1}, u_{k, k})$ by definition. It follows that $b\le u_{k, 1}< u_{k, k}\le \sigma_{\mathbf x}(b)$. 
To identify the subtrees above $b$, let us first introduce the post-$b$ process $\mathbf z^{(b)}=(z(t))_{t\ge 0}$ as follows:
\[
z^{(b)}(t)=x(b+t)-x(b-),
\]
if $0\le t\le \sigma_{\mathbf x}(b)-b$, and $z^{(b)}(t)=0$ otherwise. Clearly, the lifetime of $\mathbf z$ is $\zeta(\mathbf z)=\sigma_{\mathbf x}(b)-b$. Denote by $\underline z(t)=\inf_{s\in [0, t]}z(s)$ the running infimum of $\mathbf z$ at time $t$. For $t\in [0, \zeta(\mathbf z)]$, we next define
\begin{equation}
\label{def: g-d}
g(\mathbf z, t) = \sup\big\{s\le t: z(s-) \le \underline z(t)\big\}\vee 0 \;  \text{and}\;  d(\mathbf z, t)=\inf\big\{s>t: z(s-)\le \underline z(t)\big\}\wedge \zeta(\mathbf z),
\end{equation}
with the convention $\sup\varnothing=-\infty$ and $\inf\varnothing = \infty$. 
%
For $1\le i\le k$, write $g_{i}=g(\mathbf z, u_{k, i}-b)$ and $d_{i}=d(\mathbf z, u_{k, i}-b)$. We let 
\[
\mathbf u^{-}=\big\{u_{k, 1}: 1\le i\le k, (g_{i}, d_{i})\ne (0, \zeta(\mathbf z))\big\} = \big\{u_{k, 1}: 1\le i\le k, g_{i}>0 \text{ or } d_{i}<\zeta(\mathbf z)\big\}.
\]
We define an equivalence relationship $\sim$ on $\mathbf u^{-}$ for which $u_{k, i}\sim u_{k, j}$ if and only if $(g_{i}, d_{i})=(g_{j}, d_{j})$. Let $p$ be the number of the equivalence classes and denote by $\mathbf u^{(1)}, \mathbf u^{(2)}, \dots, \mathbf u^{(p)}$ these equivalence classes, listed in the increasing order of their least elements. 
For $1\le m\le p$ and any $u_{k, j}\in \mathbf u^{(m)}$, let $\mathbf x^{(m)}$ be the portion of $\mathbf z$ running on $[g_{j}, d_{j}]$, namely,   
\[
x^{(m)}(t)= z(t+g_{j})-z(g_{j}-)=x(t+b+g_{j})-x((b+g_{j})-), \quad \text{if } t\le d_{j}-g_{j}, 
\]
and $x^{(m)}(t)=0$ otherwise. On the event $\mathbf u^{-}\ne \varnothing$, define $\mathscr T(\mathbf x; \mathbf u_{k})$ as the following ordered rooted tree:
\[
\mathscr T(\mathbf x; \mathbf u_{k}) = \{\varnothing\}\cup\bigcup_{m=1}^{p} \theta_{m}\Big(\mathscr T(\mathbf x^{(m)}; \mathbf u^{(m)})\Big).
\]
See Fig.~\ref{fig: spanning} for some examples. 

\bigskip

For $\mathbf y=(y(t))_{0\le t\le 1}\in \mathbb D([0, 1], \R)$ satisfying $y(1)=y(0)=0$, its {\it Vervaat transform}, denoted as $\Ver(\mathbf y)$, is a c\`adl\`ag function $\mathbf x=(x(t))_{0\le t\le 1}$ defined by
\begin{equation}
\label{def: Ver-op}
x(t) = y\big(\rho_{\mathbf y}+t\!\!\!\mod 1\big) - y(\rho_{\mathbf y}-)\wedge y(\rho_{\mathbf y}), \quad 0\le t<1, 
\end{equation}
where $\rho_{\mathbf y}=\inf\{t>0: y(t)\wedge y(t-) = \inf_{0\le s\le 1}y(s)\}$ is the first infimum point of $\mathbf y$. We also set $x(1)=x(1-)$. It is then clear from the definition that $x(t)\ge 0$ for all $t\in [0, 1]$. 
Also, we have $x(0)>0$ if $\mathbf y$ jumps upwards at $\rho_{\mathbf y}$. 
Let us also note that we can recover $\mathbf y$ by splitting $\mathbf x$ at $1-\rho_{\mathbf y}$. Indeed, let us set
\begin{equation}
\label{def: ver-inverse}
\Ver^{-1}(\mathbf x, \rho_{\mathbf y}):=(\tilde y(t))_{0\le t\le 1}, \quad \text{with} \quad \tilde y(t) = x\big(t+1-\rho_{\mathbf y}\!\!\!\mod 1\big) - x(1-\rho_{\mathbf y}).
\end{equation}
Then we have $\Ver^{-1}(\mathbf x, \rho_{\mathbf y})=\mathbf y$.  

Comparing \eqref{def: Ver-op} with \eqref{def: Vervaat}, we see that $\bX^{\alpha}=\Ver(\rY^{\alpha})$. Let us denote by $\bX^{\bth}=\Ver(\rY^{\bth})$ its analogue for the extremal exchangeable bridge $\rY^{\bth}$ in \eqref{def: Ytheta}. 
Let $(U_{i})_{i\ge 1}$ be a sequence of independent variables with uniform distribution in $[0, 1]$. We will show that  the above procedure of extracting a tree from a c\`adl\`ag function, when applied separately to $\bX^{\alpha}$ and $\bX^{\bth}$, will result in the spanning trees of the stable tree and ICRT. Strictly speaking, the trees $T_{k}$ and $R'_{k}$ in Section \ref{sec: outline} are labelled rather than ordered. We therefore introduce the following labelled version of $\mathscr T$; see also Fig.~\ref{fig: spanning}.  

\smallskip
\noindent
{\bf Labelled spanning trees. }
Let us denote by $\hat{\mathscr T}^{\alpha}_{k}$ the rooted  graph tree that shares the same shape as $\mathscr T(\bX^{\alpha}\,; \{U_{1}, U_{2}, \dots, U_{k}\})$: the  vertex set of $\hat{\mathscr T}^{\alpha}_{k}$ consists of the elements of the latter; $\{u, v\}$ is an edge of $\hat{\mathscr T}_{k}$ if and only if $u$ is the parent of $v$ or $v$ is the parent of $u$ in $\mathscr T(\bX^{\alpha}\,; \{U_{1}, U_{2}, \dots, U_{k}\})$; root the tree at the vertex $\varnothing$. 
If $\hat{\mathscr T}^{\alpha}_{k}$ has fewer than $k$ leaves (root excluded), use the symbol $\partial$ to denote a cemetery state and set $\mathscr T^{\mathrm{lab}}(\bX^{\alpha}\,; \{U_{1}, U_{2}, \dots, U_{k}\}) =\partial$. Otherwise, assign a uniform labelling of $1, 2, 3, \dots, k$ to the $k$ leaves of $\hat{\mathscr T}^{\alpha}_{k}$. 
Attach a leaf labelled as 0 to the root and make that leaf to be the new root of the tree, so that the tree is always rooted at a leaf. 
Remove any vertex of degree 2 by merging the two edges adjacent to the vertex. 
For each branch point $b$, let $i(b)<j(b)$ be the pair of the least leaf labels so that $b$ is the most recent common ancestor of Leaf $i(b)$ and Leaf $j(b)$. Then order the branch points according to the lexicographic order on $\N^{2}$ and label them as $b_{1}, b_{2}, b_{3}$ and etc. Observe that this corresponds to the labelling rules in Fig.~\ref{fig}. 
Denote the resulting labelled tree as $\mathscr T^{\mathrm{lab}}(\bX^{\alpha}\,; \{U_{1}, U_{2}, \dots, U_{k}\})$. 
Define $\mathscr T^{\mathrm{lab}}(\bX^{\bth}\,; \{U_{1}, U_{2}, \dots, U_{k}\})$ in a similar way. 

\smallskip
The main results of this section are the following ones, whose proofs are found respectively in Section~\ref{sec: span-stable} and Section \ref{sec: spnn-icrt}.


\begin{prop}
\label{prop: sp-stable}
For each $k\ge 1$,  we have
\[
(T_{k})_{k\ge 1} \eqd \big(\mathscr T^{\mathrm{lab}}(\bX^{\alpha}\,; \{U_{1}, U_{2}, \dots, U_{k}\})\big)_{k\ge 1}\,. 
\]
\end{prop}

\begin{prop}
\label{prop: sp-icrt}
For each $k\ge 1$, we have
\[
(R'_{k})_{k\ge 1} \eqd \big(\mathscr T^{\mathrm{lab}}(\bX^{\bth}\,; \{U_{1}, U_{2}, \dots, U_{k}\})\big)_{k\ge 1}\,. 
\]
\end{prop}

\subsection{Proof of the main theorems}
\label{sec: main-pf}

Before proceeding to the proof of Propositions \ref{prop: sp-stable} and \ref{prop: sp-icrt}, let us first explain how they will lead to Theorems \ref{thm: main'} and \ref{thm: main}.

\begin{proof}[Proof of Theorem \ref{thm: main'}] 
Kallenberg's classic result (see \eqref{eq: canon} below) 
implies that the stable bridge $\rY^{\alpha}$ is a mixing of $\rY^{\bth}$. More precisely, for any measurable and positive functional $H$ of the Skorokhod space $\mathbb D([0, 1],\R)$, we have
\[
\mathbb E[H(\rY^{\alpha})] = \int_{\boldsymbol{\Theta}} \mathbb P(\boldsymbol\Delta^{\downarrow}_{\alpha}\in d\bth)\mathbb E[H(\rY^{\bth})],
\]
where $\boldsymbol\Delta^{\downarrow}_{\alpha}$ corresponds to the sequence of jumps of $\rY^{\alpha}$ ranked in decreasing order, and $\boldsymbol{\Theta}$ is the parameter space defined in \eqref{def: bTheta}. We can replace in above $\rY^{\alpha}$ by $\bX^{\alpha}$ and $\rY^{\bth}$ by $\bX^{\bth}$, as the Vervaat transformation is measurable. 
Applying this to $\mathscr T$ and its labelled version $\mathscr T^{\mathrm{lab}}$ (here we tacitly assume the sequence $(U_{k})_{k\ge 1}$ and the randomness used for leaf labelling are all defined on the same probability space), together  with Propositions \ref{prop: sp-stable} and \ref{prop: sp-icrt}, we deduce that for a measurable functional $G: \mathbf T^{\infty}_{\mathrm{discrete}}\to \R_{+}$, 
\[
\mathbb E\big[G\big((T_{k})_{k\ge 1}\big)\big] =   \int_{\boldsymbol{\Theta}} \mathbb P(\boldsymbol\Delta^{\downarrow}_{\alpha}\in d\bth)\mathbf E^{\bth}\big[G\big((R'_{k})_{k\ge 1}\big)\big],
\]
which is the statement of Theorem \ref{thm: main'}. 
\end{proof} 

\begin{proof}[Proof of Theorem \ref{thm: main}]
According to Lemma~\ref{lem: normalisation}, Propositions~\ref{prop: stable-recov} and~\ref{prop: icrt-recov}, we can find a measurable function $\mathscr S: \mathbf T_{\mathrm{discret}}^{\infty}\to \mathbb T$ so that 
\[
(\cT_{\alpha}, \alpha^{1-1/\alpha}\,d_{\alpha}, \mu_{\alpha})=\mathscr S\big(\boldsymbol{\Delta}^{\downarrow}_{\alpha}, (T_{k})_{k\ge 1}\big) \text{ under } \mathbb P \ \text{ and } \ (\cT, d, \mu)=\mathscr S\big(\bth, (R'_{k})_{k\ge 1}\big) \text{ under } \bP^{\bth}.
\]
Applying Theorem \ref{thm: main'} to $G=F\circ \mathscr S$ concludes the proof.
\end{proof}

\subsection{Spanning trees of the stable tree}
\label{sec: span-stable}

We prove Proposition~\ref{prop: sp-stable} here by comparing our definition of $\mathscr T(\mathbf x, \mathbf u_{k})$ with the one in Section 3.2.1 of \cite{DuLG02}. 
Fix $0<u_{1}<u_{2}<\cdots<u_{k}<1$ and recall the height process $\mathrm H$ for $\bX^{\alpha}$. 
We briefly recall from \cite{DuLG02} the following definition of an ordered rooted tree with $k$ leaves as a function of $\mathrm H$ and $\mathbf u=\{u_{1}, u_{2}, \dots, u_{k}\}$. 

\smallskip
\noindent
{\bf Defining the spanning trees from the height process.} If $k=1$, set $\tilde{\mathscr T}(\mathrm H; \mathbf u) =\{\varnothing\}$. For $k\ge 2$,  let 
\[
b'=\inf\Big\{t<u_{1}: \inf_{s\in[t, u_{1}]}H_{s}\ge  \inf_{s\in [u_{1}, u_{k}]}H_{s}\Big\}.
\]
Note that a.s.~we have $b'< u_{1}$ and $H_{b'}=\inf_{s\in [u_{1}, u_{k}]}H_{s}$ as $\mathrm H$ has continuous sample paths. 
For each $u_{i}$, let $(\tilde g_{i}, \tilde d_{i})$ be the excursion interval of $\mathrm H$ above the level $H_{b'}$, namely, 
\[
\tilde g_{i}=\sup\{s<u_{i}: H_{s}=H_{b'}\} \quad\text{and}\quad \tilde d_{i}=\inf\{s>u_{i}: H_{s}=H_{b'}\}. 
\]
Say $u_{i}$ and $u_{j}$ are equivalent if they share the same excursion interval: $(\tilde g_{i}, \tilde d_{i})=(\tilde g_{j}, \tilde d_{j})$, and denote by $\tilde{\mathbf u}^{(1)}, \tilde{\mathbf u}^{(2)}, \dots, \tilde{\mathbf u}^{(\tilde p)}$ the equivalence classes of this equivalence relation. For $1\le m\le \tilde p$ and supposing $u_{j}\in\tilde{\mathbf u}^{(m)}$, let $\mathrm H^{(m)}=(H^{(m)}_{s})_{s\ge 0}$ be defined as 
\[
H^{(m)}_{s}=H_{s+\tilde g_{j}}-H_{\tilde d_{j}}, \quad \text{if } 0\le s\le \tilde d_{j}-\tilde g_{j},
\]
and $H^{(m)}_{s}=0$ otherwise. 
Let $\tilde{\mathscr T}(\mathrm H; \mathbf u)$ be the ordered rooted tree defined by
\[
\tilde{\mathscr T}(\mathrm H; \mathbf u) = \{\varnothing\}\cup\bigcup_{m=1}^{\tilde p} \theta_{m}\Big(\tilde{\mathscr T}(\mathrm H^{(m)}; \tilde{\mathbf u}^{(m)})\Big).
\]

Applying the definition \eqref{def: mrca} to $\bX^{\alpha}$, we have 
\[
b=\tau\Big(\bX^{\alpha}, u_{1}, \inf_{s\in [u_{1}, u_{k}]}X^{\alpha}_{s}\Big) = \inf\Big\{t<u_{1}: \inf_{s\in [t, u_{1}]}X^{\alpha}_{s}\ge \inf_{s\in [u_{1}, u_{k}]}X^{\alpha}_{s}\Big\}.
\]
Note that almost surely $\Delta X^{\alpha}_{u_{1}}=0$ and $\inf_{s\in [u_{1}, u_{k}]}X^{\alpha}_{s}<X^{\alpha}_{u_{1}}$. Therefore, $b< u_{1}$ a.s. For $t\ge 0$, let us define the post-$t$ process $\mathrm Z^{(t)}$ as $Z^{(t)}_{s}=X^{\alpha}_{s+t}-X^{\alpha}_{t-}$, for $0\le s\le \sigma_{\bX^{\alpha}}(t)-t$, and $Z^{(t)}_{s}=0$ for $s>\sigma_{\bX^{\alpha}(t)}-t$.  We also denote $\underline Z^{(t)}_{ \,s}=\inf_{u\in [0, s]}Z^{(t)}_{u}$. 
We start with the  following observation. 

\begin{lem}
\label{lem: same-root}
We have $b=b'$ almost surely. Moreover, the connected components of $\{t\in [0, \sigma_{\bX^{\alpha}}(b)-b]: Z^{(b)}_{t}>\underline Z^{(b)}_{\,t}\}$ coincide with the connected components of $\{t\in [0, \sigma_{\bX^{\alpha}}(b)-b]: H_{t+b}> H_{b}\}$. 
\end{lem}

\begin{proof}
Let $\widehat X^{\alpha}_{s}=X^{\alpha}_{u_{1}}-X^{\alpha}_{(u_{1}-s)-}$, for $s\in [0, u_{1}]$. We note that 
\[
u_{1}-b=\sup\Big\{t\ge 0: \sup_{s\in [0, t]}\widehat X^{\alpha}_{s} < X^{\alpha}_{u_{1}}-\!\!\inf_{s\in [u_{1}, u_{k}]}\!\!X^{\alpha}_{s}\Big\}
\overset{\text{a.s.}}{=}\inf\Big\{t>0: \sup_{s\in [0, t]}\widehat X^{\alpha}_{s}\ge X^{\alpha}_{u_{1}}-\!\!\inf_{s\in [u_{1}, u_{k}]}\!\!X^{\alpha}_{s}\Big\}.
\]
Properties of stable processes  imply that $t\mapsto\sup_{s\in [0, t]}\widehat X^{\alpha}_{s}$ only increases by jumps. Therefore,  $u_{1}-b$ is a.s.~a jump time of $\widehat{\bX}^{\alpha}$, i.e.~$b$ is a jump time of $\bX^{\alpha}$. We claim that $b'$ is also a jump time of $\bX^{\alpha}$. Indeed, 
let $t_{0}=\inf\{t>u_{1}: H_{t}=m(\mathrm H, u_{1}, u_{k})\}$; then a.s.~$t_{0}\in (u_{1}, u_{k})$ and is therefore a local minimum point of $\mathrm H$. 
From the encoding \eqref{def: H-encoding} we can readily check  that its projection onto the tree $p(t_{0})$ is a branch point, and from its definition we have $b'=\min p^{-1}(\{p(t_{0})\})$. It is then a well known property of stable trees that $b'$ is a jump time of $\bX^{\alpha}$ (see for instance \cite{DuLeG05}, Theorem 4.6). 
Suppose that $\tau$ is a jump time of $\bX^{\alpha}$; recall the post-$\tau$ process 
$\mathrm Z^{(\tau)}$.
Now let us show that 
\begin{align*}
&\text{(S) \quad \  $H_{t}\ge H_{\tau}$ for all $\tau\le t\le  \sigma_{\bX^{\alpha}}(\tau)$; moreover 
the connected components of $\{t\in [0,\sigma_{\bX^{\alpha}}(\tau)-\tau]: $}\\
&\qquad \quad \text{$Z^{(\tau)}_{t}>\underline Z^{(\tau)}_{\,t}\}$ coincide with the connected components of $\{t\in [0,\sigma_{\bX^{\alpha}}(\tau)-\tau]: H_{\tau+t}> H_{\tau}\}$.}
\end{align*}
 Appealing to the excursion theory and the scaling property, we only need to show this for the stable process $\rY^{\alpha}$. However, for $\rY^{\alpha}$, the definition \eqref{def: H} of the height process implies that the excursion intervals of $\mathrm Y^{\alpha}$ above its running infimum $\mathrm I$ coincide with those of $\mathrm H$ away from $0$. 
 Fix $\epsilon>0$ and let $\tau^{1}_{\epsilon}$ be the first moment $t$ such that $\Delta Y^{\alpha}_{t}\ge\epsilon$. Strong Markov property and the previous arguments imply that (S) holds true for $\tau^{1}_{\epsilon}$. Repeatedly apply this arguments to the successive jump moments and then let $\epsilon\to 0$. This leads to the desired result. 
In particular, (S) implies that almost surely
 \begin{equation}
 \label{eq: ht-x}
 H_{t}\ge H_{\tau}, \quad \text{for all } t\in \big[\tau, \sigma_{\bX^{\alpha}}(\tau)\big]
 \end{equation}
 for every jump time $\tau$ of $\bX^{\alpha}$. We have already seen that the definition of $b$ ensures that $b\le u_{1}<u_{k}\le \sigma_{\bX^{\alpha}}(b)$. 
If $b'< b$, then applying \eqref{eq: ht-x} respectively to $b'$ and $b$, we find that $H_{b'}\le \inf_{u\in [b', b]}H_{u} \le H_{b}\le \inf_{u\in [u_{1}, u_{k}]}H_{u}=H_{b'}$, which implies that $H_{b}=\inf_{u\in [b', b]}=H_{b}$. Let us briefly argue that this occurs with null probability. Since $\tau_{\bX^{\alpha}}(b)$ is a stopping time, we have $\inf_{u\in [\tau, \tau'+\epsilon]} H_{u}<H_{\tau}$ a.s.~for all $\epsilon>0$ (Lemma 1.4.5 of \cite{DuLG02}), where $\tau=\tau_{\bX^{\alpha}}(b)$. Combined with the time-reversal property of $\mathrm H$ (Corollary 3.1.6 of \cite{DuLG02}), we see that $\inf_{u\in [b-\epsilon, b]}H_{u}<H_{b}$. 
Hence, we must have $b\le b'$. To show the other side, let $t_{1}\in [u_{1}, u_{k}]$ be such that $X^{\alpha}_{t_{1}-}=\inf_{s\in [u_{1}, u_{k}]}X^{\alpha}_{s}$. If $s\in [0, t_{1})$ satisfies $X^{\alpha}_{s-}< \inf_{u\in [s, t_{1}]}X^{\alpha}_{u}$, then we must have $s\le b$, since $\inf_{u\in [b, t_{1}]}X^{\alpha}_{u}\ge X^{\alpha}_{t_{1}-}$ by the choice of $b$. 
Now take $r\in [b, t_{1}]$; the previous arguments imply that 
\[
\Big\{s\in [0, t_{1}): X^{\alpha}_{s-}< \inf_{u\in [s, t_{1}]}X^{\alpha}_{u}\Big\} \subseteq \Big\{s\in [0, r): X^{\alpha}_{s-}< \inf_{u\in [s, r]}X^{\alpha}_{u}\Big\}.
\]
It then follows from \eqref{def: H} that $H_{t_{1}}\le \min_{r\in [b, t_{1}]}H_{r}$. Compared with the definition of $b'$, this suggests that $b'\le b$. We conclude with $b=b'$. The second part of the lemma follows from (S). 
\end{proof}

\begin{proof}[Proof of Proposition~\ref{prop: sp-stable}]
It suffices to show that 
\[
\big(\tilde{\mathscr T}(\mathrm H; \{U_{1}, U_{2}, \dots, U_{k}\})\big)_{k\ge 1}\eqd \big(\mathscr T(\bX^{\alpha}; \{U_{1}, U_{2}, \dots, U_{k}\})\big)_{k\ge 1}\,.
\]
since $T_{k}$ is obtained from $\tilde{\mathscr T}(\mathrm H; \{U_{1}, U_{2}, \dots, U_{k}\})$. 
Moreover, it suffices to prove the above identity in distribution for each $k$, 
since $\mathscr T(\bX^{\alpha}; \{U_{1},\dots, U_{k}\})$ is a subtree of $\mathscr T(\bX^{\alpha}; \{U_{1},\dots, U_{k+1}\})$. 
For each realisation of $\{U_{1},\dots, U_{k}\}$, Lemma~\ref{lem: same-root} says that the root degree in both trees are the same and the intervals that will be used to build the subtrees above also coincide. This is enough to conclude thanks to the recursive nature of both definitions.
\end{proof}

\section{Spanning trees of the ICRT}
\label{sec: spnn-icrt}

In this section, we prove Proposition~\ref{prop: sp-icrt} using weak convergence arguments. We introduce in Section~\ref{sec: intro-ptree} the counterpart of the ICRT in the discrete world: the model of $\bp$-trees. In Section~\ref{sec: LIFO-ptree} we describe an encoding of $\bp$-trees. When plugging these coding processes into the function $\mathscr T(\cdot\,; \{U_{1}, U_{2}, \dots, U_{k}\})$, we obtain the spanning trees of $\bp$-trees. 
 Moreover, we will see that these coding processes converge to the extremal exchangeable process $\bX^{\bth}$ in a suitable regime, whilst the spanning trees of the $\bp$-trees will converge to those of the ICRT. We show in Section~\ref{sec: cv-icrt-sp} the function $\mathscr T$ also converges alongside the coding processes, which then allows us to  conclude the proof of Proposition~\ref{prop: sp-icrt}. 

\subsection{Preliminaries on $\bp$-trees}
\label{sec: intro-ptree}

Let $\bp_{n} = (p_{n}(i))_{1\le i\le n}$ be a probability measure on $[n]:=\{1, 2, \dots, n\}$. We further assume that $p_{n}(1)\ge p_{n}(2)\ge \cdots\ge p_{n}(n)>0$. We view a rooted tree as a family tree: the root is the common ancestor, its neighbours are the first generation, and so on. 
Denote by $\mathbf T_{n}$ the set of all labelled rooted trees with the set of vertex labels given by $[n]$. 
Cayley's multinomial formula (\cite{Pi99}) says that the following is a probability measure on $\mathbf T_{n}$:  
\begin{equation}
\label{def: p-tree-dist}
\pi^{\bp_{n}}(t):=\prod_{i\in [n]}p^{\kappa_{i}}_{n}(i), \quad t\in \mathbf T_{n}, 
\end{equation}
where $\kappa_{i}=\kappa_{i}(t)$ is the number of children of the vertex $i$ in $t$. A random tree is called a {\it$\bp_{n}$-tree} if its law is $\pi^{\bp_{n}}$. 
We are interested in the large-size limit of these trees. More precisely, the relevant asymptotic regime is as follows: suppose that there exists $\bth=(\theta_{i})_{i\ge 1}\in \boldsymbol\Theta$ so that 
\begin{equation}
\label{def: p-asym}
\frac{p_{n}(i)}{\sigma_{n}} \xrightarrow{n\to\infty} \frac{\theta_{i}}{\|\bth\|}, \quad i\ge 1, \quad\text{with}\quad \sigma_{n}:=\Big(\sum_{1\le i\le n}p^{2}_{n}(i)\Big)^{1/2}\to 0.
\end{equation}
Let $T_{n}$ be a $\bp_{n}$-tree. We turn it into a measured metric space by equipping it with the graph distance $\dgr$ and the probability measure $\bp_{n}$ on its vertex set. 
Camarri and Pitman \cite{Pi00} show that (taking into account the scaling relation \eqref{id: icrt-scaling})
\begin{equation}
\label{cv: p-tree}
\big(T_{n}, \sigma_{n}\dgr, \bp_{n}\big) \xrightarrow[n\to\infty]{(d)} (\cT, \|\bth\|\cdot d, \mu) \ \text{under } \bP^{\bth},
\end{equation}
with respect to the Gromov--Prokhorov topology. The original result in \cite{Pi00} was stated in terms of the convergence of spanning trees, which will also be useful later. More precisely, for each $n$, let $(\eta^{(n)}_{k})_{k\ge 1}$ be a sequence of independent variables with common distribution $\bp_{n}$ and denote by $\hat{R}^{n}_{k}$ the smallest subtree of $T_{n}$ containing the root and the vertices  $\eta^{(n)}_{1}, \eta^{(n)}_{2}, \dots, \eta^{(n)}_{k}$. 
If $\eta^{(n)}_{1}, \eta^{(n)}_{2}, \dots, \eta^{(n)}_{k}$ are not distinct leaves in $\hat{R}^{n}_{k}$, set $R^{n}_{k}=\partial$. Otherwise, relabel  $\eta^{(n)}_{1}, \eta^{(n)}_{2}, \dots, \eta^{(n)}_{k}$ uniformly from $1$ to $k$ and remove any vertex of degree 2 in $\hat{R}^{n}_{k}$. Relabel the root as 0. For each branch point $b$, let $i(b)<j(b)$ be the pair of the least leaf labels so that $b$ is the most recent common ancestor of Leaf $i(b)$ and Leaf $j(b)$. Then order the branch points according to the lexicographic order on $\N^{2}$ and label them as $b_{1}, b_{2}, b_{3}$ and etc. 
Call the resulting tree $R^{n}_{k}$. 
Observe that $R'_{k}$, the graph tree obtained from the spanning tree $\cR'_{k}$ of $\cT$, is labelled in the same way. Clearly, the set of all graph trees with $k$ leaves labelled from 0 to $k$, root at 0, no vertices of degree 2, and branch points labelled as $(b_{i})_{i\ge 1}$ contains finite elements. 
Then Camarri and Pitman \cite{Pi00} show that for each $k\ge 1$, 
\begin{equation}
\label{cv: p-sp-tree}
\lim_{n\to\infty}\mathbb P(R^{n}_{k}=\partial) =0, \ \text{and the law of $R^{n}_{k}$ is equal to that of $R'_{k}$ for $n$ sufficiently large. } 
\end{equation}

\subsection{A LIFO queue construction of $\bp$-trees}
\label{sec: LIFO-ptree}

Let $(\chi_{i})_{i\ge 1}$ be a sequence of independent uniform variables on $[0, 1]$. For each $n\ge 1$, consider an exchangeable process $\rY^{n}=(Y^{n}_{t})_{0\le t\le 1}$ defined as follows:
\begin{equation}
\label{def: Yn}
Y^{n}_{t}=-t + \sum_{1\le i\le n }p_{n}(i)\mathbf 1_{\{\chi_{i}\le t\}}=\sum_{1\le i\le n}p_{n}(i)\big(\mathbf 1_{\{\chi_{i}\le t\}}-t\big), \quad 0\le t\le 1.
\end{equation}
Performing the Vervaat transformation \eqref{def: Ver-op} on $\rY^{n}$ results in an excursion-like process $\bX^{n}=\Ver(\rY^{n})$. 
Note that the jump times of $\bX^{n}$ are $\chi'_{i}=\chi_{i}-\rho_{n}\mod 1$, where $\rho_{n}$ is the first infimum point of $\rY^{n}$, $1\le i\le n$. 
Recall from Section \ref{sec: coding} the mapping $\mathscr T$, which extracts an ordered rooted tree from an excursion-like c\`adl\`ag function.  Denote by $T^{\mathrm{ord}}_{n}=\mathscr T(\bX^{n}; \{\chi'_{1}, \chi'_{2}, \dots, \chi'_{n}\})$. 
We will show below that $T^{\mathrm{ord}}_{n}$ is an ordered version of the $\bp_{n}$-tree. Before launching the proof, let us point out that this statement is implied in the Remark in Section 3.2 of \cite{AMP04}. We provide here a proof that highlights the connection between $\bp$-trees and Bienaym\'e trees. 
%
%
%

\paragraph{A LIFO queue construction for random trees.}
As a first step in identifying the distribution of $T^{\mathrm{ord}}_{n}$, we explain here an alternative construction of the tree. Imagine a queuing system with a single server and $n$ customers $1, 2, \cdots, n$. Customer $i$ enters the queue at time $\chi'_{i}$ and requires the attention of the server for an amount $p_{n}(i)$ of service time, $1\le i\le n$. The server operates under a Last-In-First-Out (LIFO) rule. That is, when a new customer arrives, the server immediately interrupts the current service and serves the new arrival. Only after the new arrival leaves the queue does the server come back to the last customer in the queue. It is then not difficult to check that $X^{n}_{t}$ is the amount of unfulfilled service time for the customers in the queue (i.e.~load of the server) at time $t$. 

Now introduce a genealogy on the customers by declaring the first arriving customer as the root; moreover Customer $j$ is a child of Customer $i$ if and only if the former interrupts the service of the latter. Note that Customer $i$ leaves the queue at time $\sigma_{\bX^{n}}(\chi'_{i})$, namely the first moment when the load of the server falls back to the level prior to its arrival. 
Thus, the descendants of Customer $i$ are those who arrive between $\chi'_{i}$ and $\sigma_{\bX^{n}}(\chi'_{i})$. 
We further assume that these descendants are ranked in their arrival orders.  
Let us note that $\bX^{n}$ only increases by jumps; thus a branch point of $\mathscr T(\bX^{n}; \{\chi'_{1}, \chi'_{2}, \dots, \chi'_{n}\})$ must correspond to a jump of $\bX^{n}$. 
It can be checked from its definition in Section \ref{sec: coding} that 
the genealogy on the $n$ jumps of $\bX^{n}$ obtained from the LIFO-queue is the same as the one given by $T^{\mathrm{ord}}_{n}:=\mathscr T(\bX^{n}; \{\chi'_{1}, \chi'_{2}, \dots, \chi'_{n}\})$. 
For later discussion, it will be important to retain the information on the service times. Therefore, we label the vertex corresponding to the jump at time $\chi'_{i}$ as $i$ and assign the mark $p_{n}(i)$ to it. The obtained labelled rooted tree is denoted as $T^{\mathrm{lab}}_{n}$ and we refer to $(T^{\mathrm{lab}}_{n}, (p_{n}(i))_{i\in [n]})$ as the marked labelled tree obtained from the LIFO-queue construction. In what follows, we show that $T^{\mathrm{lab}}_{n}$ is a $\bp_{n}$-tree. 
\begin{prop}
\label{prop: p-tree}
$T^{\mathrm{lab}}_{n}$ has the distribution $\pi^{\bp_{n}}$ defined in \eqref{def: p-tree-dist}. 
\end{prop}

\paragraph{Trees encoded by compound Poisson processes.}
The above LIFO-queue construction was initially introduced for excursions of compound Poisson processes in \cite{LGLJ98}. The trees obtained in this way have a remarkably  simple distribution thanks to the Markov property of the underlying process. Let us briefly explain this. 
Instead of the exchangeable process $\rY^{n}$, consider this time a process on $\R_{+}$ which is defined using the following random variables. Let $(\Delta_{i})_{i\ge 1}$ be a sequence of i.i.d.~positive random variables with a common distribution $f(x)dx$, where $f$ is a continuous probability density function with support on $[0, 1]$. In particular, we have $\mathbb E[\Delta_{1}]\le 1$. Let $0<E_{1}<E_{2}<\dots$ be the jump times in a Poisson process of unit rate. Set
\[
S_{t} = -t+ \sum_{i\ge 1} \Delta_{i}\mathbf 1_{\{E_{i}\le t\}}, \quad t\ge 0.
\]
Denote 
\[
\sigma_{\mathrm S}(E_{1}) = \inf\{t>E_{1}: S_{t} \le S_{E_{1}-}\}.
\]
Then $(E_{1}, \sigma_{\mathrm S}(E_{1}))$ is the first excursion interval of $(S_{t})_{t\ge 0}$ away from its running infimum (note that each such excursion must start with a jump). Denote by $\mathcal E=(\mathcal E_{t})_{0\le t\le \zeta}$ this excursion:
\[
\mathcal E_{t} = S_{E_{1}+t}-S_{E_{1}-}, \quad 0\le t\le\zeta:= \sigma_{\mathrm S}(E_{1})-E_{1}.
\]
Denote by $E'_{1}=0, E'_{2}=E_{2}-E_{1}, \cdots, E'_{p}=E_{p}-E_{1}$ the sequence of jump times of $\mathcal E$.  
Let $\tau=\mathscr T(\mathcal E; \{E'_{1}, E'_{2}, \dots, E'_{p}\})$, and for each $u\in \tau$, set $m_{u}=\Delta_{i}$ if $u$ corresponds to the customer arriving at $E'_{i}$. 
Le Gall and Le Jan have shown in \cite{LGLJ98} that
\begin{quote}
{\it $m_{\varnothing}$ is distributed as $\Delta_{1}$; conditional on it being $x$, the number of its children is a Poisson variable with mean $x$; the marked subtrees above the root are i.i.d.~with the same distribution as $(\tau, (m_{u})_{u\in \tau})$.}
\end{quote}
Denote by $\mathbf m^{\downarrow}$ the list of marks in $\tau$ ranked in a decreasing order. Assign the labels from $[n]$ to the vertices so that the vertex with label 1 has the largest mark, the one with label 2 has the second largest mark, and etc. 
Let $\tau^{\mathrm{lab}}$ stand for this labelled version of $\tau$ (vertex ordering is ignored). 
For any $t\in \mathbf T_{n}$ and any Borel set $B\subset \{(x_{1}, x_{2}, \dots, x_{n})\in [0, 1]^{n}: x_{1}\ge x_{2}\ge \cdots \ge x_{n}\}$, we have
\begin{equation}
\label{eq: mix-p-tree}
\mathbb P\Big(\tau^{\mathrm{lab}}=t; \mathbf m^{\downarrow}\in B\Big) = 
\int_{B} e^{-\sum_{i}x_{i}}\prod_{i\in [n]}f(x_{i}) x_{i}^{\kappa_{i}(t)}dx_{1}dx_{2}\cdots dx_{n}.
\end{equation}
{\it Proof of \eqref{eq: mix-p-tree}}: for $t\in \mathbf T_{n}$, there are $\prod_{i}(k_{i}(t))!$ ways of ordering its vertices. Combining this with the aforementioned result of Le Gall and Le Jan \cite{LGLJ98}, we deduce the formula.  

As a consequence of \eqref{eq: mix-p-tree}, we note that conditioned on $\mathbf m^{\downarrow}=\bp_{n}$, $\tau^{\mathrm{lab}}$ is distributed as a $\bp_{n}$-tree.

\begin{figure}[tp]
\centering
\includegraphics[height = 5cm]{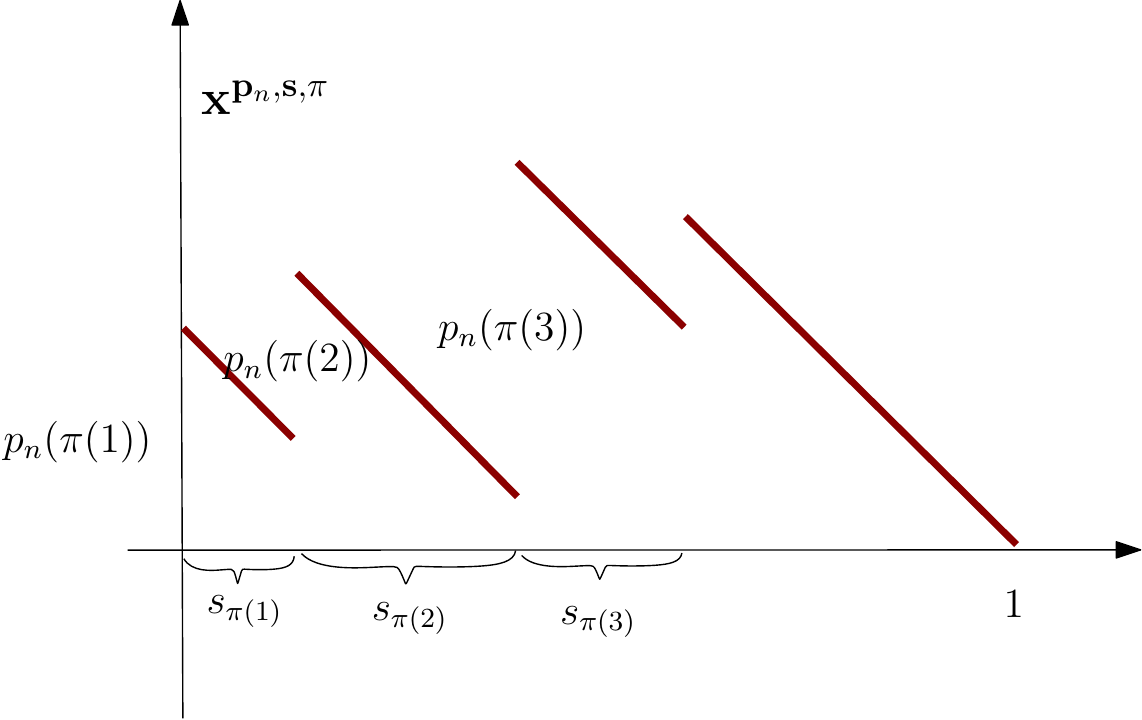}
\caption{\label{fig: path} An example of an admissible $\mathbf x^{\mathbf p_{n}, \mathbf s, \pi}$. Successive jump sizes are $p_{n}(\pi(1)), p_{n}(\pi(2)), p_{n}(\pi(3)), \cdots$, and the gaps between the jumps are given by $s_{\pi(1)}, s_{\pi(2)}, s_{3(\pi(3)}, \cdots$.}
\end{figure}

\paragraph{Excursion of $\mathrm S$ conditioned on its jumps.}
Recall $\mathcal E$ has jumps at $0=E'_{1}< E'_{2} < \cdots < E'_{p}$ and 
recall that $\mathbf m^{\downarrow}$ is also the sequence of jump sizes in $\mathcal E$ listed in decreasing order. Set $E'_{p+1}=1$. 
We observe that $\mathcal E$ is characterised by its jump sizes at $(E'_{i})_{1\le i\le p}$ and $(E'_{i+1}-E'_{i})_{1\le i\le p}$. 
Given a permutation $\pi$ of $[n]$ and a sequence $\mathbf s=(s_{1}, s_{2}, \cdots, s_{n-1}, s_{n})\in [0, 1]^{n}$ satisfying $\sum_{1\le i\le n}s_{i}=1$,  
we define a path $\mathbf x^{\bp_{n},\mathbf s, \pi}\in \mathbb D([0, 1], \R)$ as follows: for $0\le t<1$, let
\[
x^{\bp_{n},\mathbf s, \pi}(t)= \sum_{j=1}^{i(t)}p_{n}(\pi(j))-t, \quad \text{where } i(t)=\min\Big\{i\ge 1:  \sum_{j=1}^{i}s_{\pi(j)}>t\Big\}, 
\]
and $x(1)=x(1-)=0$. See Fig.~\ref{fig: path} for an example. 
We say $\mathbf x^{\bp_{n},\mathbf s, \pi}$ is {\it admissible} if it only takes non negative values. 
Suppose that $\mathbf x^{\bp_{n},\mathbf s, \pi}$ is admissible. From the memoryless properties of exponential variables, we deduce that the law of $\mathcal E$ has the following density at $\mathbf x^{\bp_{n},\mathbf s, \pi}$:
\begin{equation}
\label{eq: law-E}
\prod_{i\in [n]}f\big(p_{n}(i)\big) e^{-s_{i}} = e^{-1}\prod_{i\in [n]}f\big(p_{n}(i)\big)
\end{equation}
It follows that conditional on $\mathbf m^{\downarrow}=\bp_{n}$, $\mathcal E$ is uniformly distributed on the set of admissible $\mathbf x^{\bp_{n},\mathbf s, \pi}$. 

On the other hand, let us denote by $\mathbf D_{n}=\{\mathbf s=(s_{i})_{1\le i\le n}\in [0, 1]^{n}:  \sum_{1\le i\le n}s_{i}=1\}$ the $n$-dimensional simplex. We write  
\[
\mathscr D=\Big\{\mathbf s\in \mathbf D_{n}:  \exists\,\pi:[n]\to [n] ,  \exists\, t\in (0, 1):  x^{\bp_{n}, \mathbf s, \pi}(t)=0 \text{ and } x^{\bp_{n}, \mathbf s, \pi}(s)\ge 0, \forall\, s\in [0, 1]\Big\}. 
\]
Note that this is a subset of 
\[
\bigcup_{1\le i\le n-1}\bigcup_{\pi}\Big\{\mathbf s\in \mathbf D_{n}: \sum_{j=1}^{i}p_{n}(\pi(j))=\sum_{j=1}^{i}s_{\pi(j)}  \Big\}
\]
where the second union is over all the permutations $\pi$ of $[n]$. As a consequence, $\mathscr D$ has null measure under the uniform distribution on $\mathbf D_{n}$. Meanwhile, for $\mathbf s\in \mathbf D_{n}\setminus\mathscr D$, among its  $n$ cyclic permutations, there is precisely one which makes $\mathbf x^{\bp_{n}, \mathbf s, \pi}$ admissible. We'll use this to find the distribution of $\bX^{n}$. To that end, we observe that $\rY^{n}$ has jumps at $\chi_{1}, \chi_{2}, \dots, \chi_{n}$, whose joint distribution determines the law of $\rY^{n}$. Let $\chi_{n, 1}<\chi_{n, 2}< \cdots<\chi_{n, n}$ be the order statistics of $(\chi_{i})_{1\le i\le n}$. Set $\mathbf r=(r_{i})_{1\le i\le n}$ with 
\[
r_{1}=\chi_{n, 1}+1-\chi_{n, n}, \ r_{2}=\chi_{n, 2}-\chi_{n, 1}, \ r_{3}=\chi_{n, 3}-\chi_{n, 2},\  \dots\  ,r_{n}=\chi_{n, n}-\chi_{n, n-1}.
\]
It is straightforward to check that $\mathbf r$ follows the uniform distribution on $\mathbf D_{n}$. Since $\bX^{n}=\mathbf x^{\bp_{n}, \mathbf r, \pi}$ for certain cyclic permutation $\pi$, it then follows from the previous arguments that $\bX^{n}$ is uniformly distributed on the set of all admissible $x^{\bp_{n}, \mathbf s, \pi}$. Compared with \eqref{eq: law-E}, this shows that $\mathcal E$ conditioned on $\mathbf m^{\downarrow}=\bp_{n}$ has the same distribution as $\bX^{n}$. 

\begin{proof}[Proof of Proposition~\ref{prop: p-tree}]
On the one hand, \eqref{eq: mix-p-tree} says that $\tau^{\mathrm{lab}}$ conditioned on $\mathbf m^{\downarrow}=\bp_{n}$ is distributed as a $\bp_{n}$-tree. On the other hand, the previous arguments show that $\mathcal E$ conditioned on $\mathbf m^{\downarrow}=\bp_{n}$ has the same distribution as $\bX^{n}$. 
Comparing this with the definition of $T^{\mathrm{lab}}_{n}$, we see that it is distributed as $\tau^{\mathrm{lab}}$ conditioned on $\mathbf m^{\downarrow}=\bp_{n}$, and therefore a $\bp_{n}$-tree. 
\end{proof}

\noindent
{\bf Remark.}
Eq.~\eqref{eq: mix-p-tree} shows that certain types of Bienaym\'e trees are mixtures of $\bp$-trees. Since stable trees and ICRT are respectively scaling limits of Bienaym\'e trees and $\bp$-trees, it is very tempting to prove Theorem \ref{thm: main} via the weak convergence arguments. However, for that to work, at the very least we need to show that $\boldsymbol\Delta^{\downarrow}$ appears as the same functional of $\cT^{\alpha}$ as $\bth$ for the ICRT $\cT$, which is not obvious. In the current approach, this is covered by Propositions~\ref{prop: stable-recov} and \ref{prop: icrt-recov}. 

\subsection{Spanning trees of the $\bp$-trees}

Recall from Section~\ref{sec: LIFO-ptree} the $\bp_{n}$-tree $T^{\mathrm{lab}}_{n}$, where the vertices are labelled using the jumps in $\bX^{n}$. Let $(\eta^{(n)}_{i})_{i\ge 1}$ be i.i.d.~variables with common law $\bp_{n}$, which we view as a distribution on the vertex set of $T^{\mathrm{lab}}_{n}$. Thanks to Proposition~\ref{prop: p-tree}, we know that the subtree of $T^{\mathrm{lab}}_{n}$ spanned by $\eta^{(n)}_{1}, \eta^{(n)}_{2}, \dots, \eta^{(n)}_{k}$ has the same law as the spanning tree $R^{n}_{k}$ of $\bp_{n}$-trees, seen in Section~\ref{sec: intro-ptree}. Abusing the notation, we denote this subtree of $T^{\mathrm{lab}}_{n}$ as $R^{n}_{k}$. 
Recall that $(U_{i})_{i\ge 1}$ is an i.i.d.~sequence of uniform variables on $(0, 1)$, independent of $\rY^{n}$. Follow the same rules as set out in the paragraph {\bf Labelled spanning trees} in Section \ref{sec: coding} to obtain a labelled version $\mathscr T^{\mathrm{lab}}\big(\bX^{n}; \{U_{1}, U_{2}, \dots, U_{k}\})$ of $\mathscr T\big(\bX^{n}; \{U_{1}, U_{2}, \dots, U_{k}\})$. 
Let us show the following. 
\begin{prop}
\label{prop: cd-sp}
Let $k\ge 1$. Assume that \eqref{def: p-asym} is true. There exists a coupling between $(\eta^{(n)}_{i})_{1\le i\le k}$ and $(U_{i})_{1\le i\le k}$ so that for $n$ sufficiently large, we have 
\[
R^{n}_{k}= \mathscr T^{\mathrm{lab}}\big(\bX^{n}; \{U_{1}, U_{2}, \dots, U_{k}\}\big).
\]
\end{prop}

\begin{figure}
\centering
\includegraphics[height = 7.5cm]{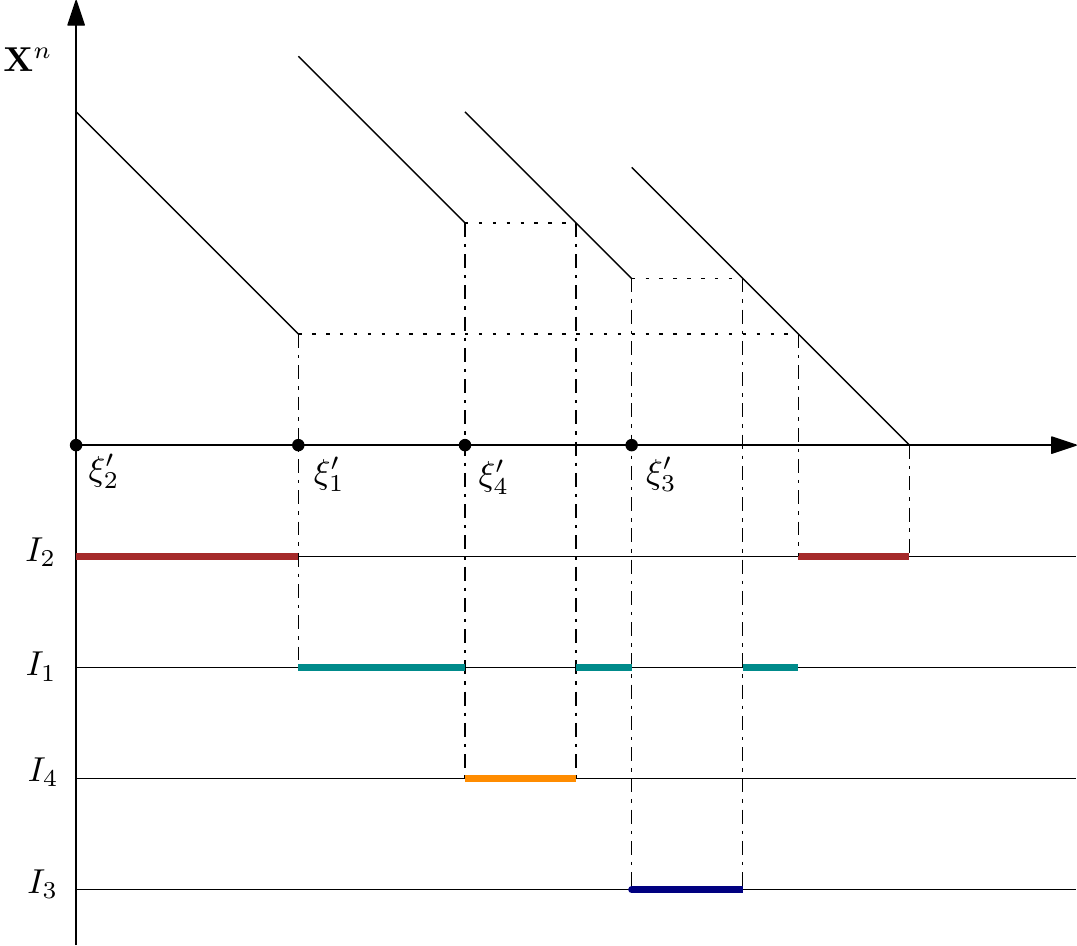}
\caption{\label{fig2}An example of the partition of $[0, 1]$ into $(I_{j})_{1\le j\le k}$.}
\end{figure}

\begin{proof}
For $t\in (0, 1)$, we define
\[
q(t) = \tau\big(\bX^{n}, t, X^{n}_{t}\big) = \inf\Big\{s<t: \inf_{u\in [s, t]}X^{n}_{u} \ge X^{n}_{t}\Big\}.
\]
Borrowing the LIFO-queue metaphor, we can say that $q(t)$ refers to the arrival time of the client that the server is serving at time $t$. 
Since $\bX^{n}$ only increases at its jump times $\chi'_{1}, \chi'_{2}, \dots, \chi'_{n}$, one can show that $\Delta X^{n}_{q(t)}>0$; thus $q(t)=\chi'_{i}$ for some $i\in [n]$. Let 
\[
\tilde R^{n}_{k}=\mathscr T^{\mathrm{lab}}\big(\bX^{n}; \{q(U_{1}), q(U_{2}), \dots, q(U_{k})\}\big).
\]
Then $\tilde R^{n}_{k}$ corresponds to the subtree of $T^{\mathrm{ord}}_{n}$ spanned by $q(U_{1}), q(U_{2}), \dots, q(U_{k})$. Let us define $I_{j}=\{t\in (0, 1): q(t) = \chi'_{j}\}$, $1\le j\le n$. Using  the fact that 
$\bX^{n}$ has the drift $-1$ and an induction on $n$, it is elementary to check (see also Fig.~\ref{fig2}) that  $(I_{j})_{1\le j\le n}$ is disjoint and $I_{j}$ has the Lebesgue measure $p_{n}(j)$. 
We now couple $(\eta^{(n)}_{i})_{i\ge 1}$ with $(U_{i})_{i\ge 1}$ by putting $\eta^{(n)}_{i}=j$ if and only if $U_{i}\in I_{j}$. It follows from this coupling that $R^{n}_{k}=\tilde R^{n}_{k}$ almost surely. 

Meanwhile, since a branch point of  $\mathscr T^{\mathrm{lab}}\big(\bX^{n}; \{U_{1}, U_{2}, \dots, U_{k}\}\big)$ must correspond to a jump time of $\bX^{n}$, we see that $\mathscr T^{\mathrm{lab}}\big(\bX^{n}; \{U_{1}, U_{2}, \dots, U_{k}\}\big)$ differs from $\tilde R^{n}_{k}$ in having at most $k$ additional leaves attached respectively to $q(U_{i})$, $1\le i\le k$. 
However, the additional leaves appear only if $q(U_{i})$ has degree $\ge 3$ (see also Fig.~\ref{fig: spanning} for an example). The latter event happens if some $I_{j}$ contains more than one element of $(U_{i})_{1\le i\le k}$. But \eqref{def: p-asym} ensures that $\max_{j\ge 1}p_{n}(j)\to 0$. The conclusion follows. 
\end{proof}

\subsection{Convergence of spanning trees}
\label{sec: cv-icrt-sp}

Throughout this section, we will use the shorthand notation $\mathbb D=\mathbb D([0, 1], \R)$. 
Our aim here is to show the following proposition. Recall that $(U_{i})_{i\ge 1}$ is an i.i.d.~sequence of uniform variables on $(0, 1)$. 

\begin{prop}
\label{prop: cv-sp}
Suppose that \eqref{def: p-asym} takes place. For each $k\ge 1$ and $n$ sufficiently large, we have
\[
\mathscr T\big(\bX^{n}; \{U_{1}, U_{2}, \dots, U_{k}\}\big)\eqd \mathscr T\big(\bX^{\bth}; \{U_{1}, U_{2}, \dots, U_{k}\}\big).
\]
\end{prop}

Proposition~\ref{prop: cv-sp} will allow us to complete the proof of Proposition~\ref{prop: sp-icrt}.

\begin{proof}[Proof of Proposition~\ref{prop: sp-icrt}]
By Proposition~\ref{prop: cd-sp} and \eqref{cv: p-sp-tree}, we deduce that after finite $n$, the distribution of $\mathscr T^{\mathrm{lab}}(\bX^{n}; \{U_{1}, U_{2}, \dots, U_{k}\})$ is identical to that of  $R'_{k}$. Comparing this with Proposition~\ref{prop: cv-sp}, we find that
\[
 \mathscr T^{\mathrm{lab}}(\bX^{\bth}; \{U_{1}, U_{2}, \dots, U_{k}\}) \eqd R'_{k}, \quad k\ge 1. 
 \]
 Since the tree $R'_{k-1}$ can be obtained from $R'_{k}$ by removing the leaf labelled $k$, and similarly the tree $\mathscr T^{\mathrm{lab}}(\bX^{\bth}; \{U_{1}, U_{2}, \dots, U_{k-1}\})$ is a deterministic function of $\mathscr T^{\mathrm{lab}}(\bX^{\bth}; \{U_{1}, U_{2}, \dots, U_{k}\})$, the conclusion follows. 
\end{proof}

All it remains now is to prove Proposition~\ref{prop: cv-sp}. To that end, we require some elementary results on Skorokhod's topology.  
These are collected in Appendix \ref{sec: topo}. 
We will also need some path properties of the exchangeable process $\rY^{\bth}$ and its Vervaat transform $\bX^{\bth}$, which are stated in Section~\ref{sec: a-exch}. 

\begin{proof}[Proof of Proposition~\ref{prop: cv-sp}]
Without much loss of generality, let us assume $\|\bth\|=1$. 
According to Theorem~\ref{thm: cv-ex}, under the assumption \eqref{def: p-asym}, the exchangeable processes $\sigma_{n}\rY^{n}$ converge in distribution to $\rY^{\bth}$ in $\mathbb D$. Combined with Lemma~\ref{lem: D0}, Theorem~\ref{thm: knight} and Proposition~\ref{prop: rY-cont}, this entails the convergence in distribution of $\sigma_{n}\bX^{n}$ to $\bX^{\bth}$.  Appealing to Skorokhod's Representation Theorem, we can assume the convergence takes place almost surely. Namely,
\[
\sigma_{n}\bX^{n}\xrightarrow[n\to\infty]{a.s.} \bX^{\bth} \quad \text{in } \mathbb D.
\]
Recall that $(U_{i})_{1\le i\le k}$ is a sequence of independent uniform points on $(0, 1)$ and that $U_{k, 1}<U_{k, 2}<\cdots<U_{k, k}$ is the order statistics of $(U_{i})_{1\le i\le k}$. 
Since $\bX^{\bth}$ is a.s.~continuous at both $U_{k, 1}$ and $U_{k, k}$, 
according to Lemma~\ref{lem: D-1}, we then have
\[
m_{n}:=\sigma_{n}\inf_{s\in [U_{k, 1}, U_{k, k}]}X^{n}_{s} \xrightarrow[n\to\infty]{a.s.} m:=  \inf_{s\in [U_{k, 1}, U_{k, k}]}X^{\bth}_{s}.
\]
Lemma~\ref{lem: X-prop} ensures that all the conditions in Lemma~\ref{lem: D1} are met by $\bX^{\bth}$, $U_{k, 1}$ and $m$, so that $b_{n}:=\tau(\sigma_{n}\bX^{n}, U_{k, 1}, m_{n})\to b:=\tau(\bX^{\bth}, U_{k, 1}, m)$ and $\sigma_{n}\Delta X^{n}_{b_{n}}\to \Delta X^{\bth}_{b}$. Let us define
\begin{equation}
\label{def: Z}
Z^{n}_{s}= X^{n}_{b_{n}+(1-b_{n})s}-X^{n}_{b_{n}} \quad \text{and} \quad Z_{s}=X^{\bth}_{b+(1-b)s}-X^{\bth}_{b}, 
\quad 0\le s\le 1.
\end{equation}
Lemma \ref{lem: D15} then ensures that 
\[
\sigma_{n}\mathrm Z^{n} := (\sigma_{n}Z^{n}_{s})_{0\le s\le 1} \xrightarrow[n\to\infty]{\mathbb D} \mathrm Z:=(Z_{s})_{0\le s\le 1}.
\]
Put 
\begin{equation}
\label{def: Utilde}
\tilde U^{n}_{k, i}= (U_{k, i}-b_{n})/(1-b_{n}), \quad \tilde U_{k, i}= (U_{k, i}-b)/(1-b), \quad 1\le i\le k. 
\end{equation}
Note we have $\tilde U^{n}_{k, i}\to \tilde U_{k, i}$ almost surely for $1\le i\le k$. Recall from \eqref{def: g-d} the definitions of $g(\mathbf z, t)$ and $d(\mathbf z, t)$. 
As $\mathrm Z$ has no negative jumps, $t\mapsto \inf_{u\in [0, t]}Z_u$ is continuous. 
Lemma \ref{lem: Z-prop} and Lemma \ref{lem: D2} combined imply that
\[
 g(\sigma_{n}\mathrm Z^{n}_{s}, \tilde U^{n}_{k, i}) \to g(\mathrm Z, \tilde U_{k, i}), \quad d(\sigma_{n}\mathrm Z^{n}_{s}, \tilde U^{n}_{k, i}) \to d(\mathrm Z, \tilde U_{k, i}), \quad 1\le i\le k.
 \]
 Note that for $i\ne j$, either $(g(\mathrm Z, \tilde U_{k, i}), d(\mathrm Z, \tilde U_{k, i}))\cap (g(\mathrm Z, \tilde U_{k, j}), d(\mathrm Z, \tilde U_{k, j}))=\varnothing$ or the two intervals are identical. This implies that the degree of the root in the spanning tree $\mathscr T(\sigma_{n}\bX^{n}; \{U_{1}, \cdots, U_{k}\})=\mathscr T(\bX^{n}; \{U_{1}, \cdots, U_{k}\})$, which corresponds to the number 
  \[
\#\Big\{ \big(g(\sigma_{n}\mathrm Z^{n}_{s}, \tilde U^{n}_{k, i}), d(\sigma_{n}\mathrm Z^{n}_{s}, \tilde U^{n}_{k, i})\big): 1\le i\le k\Big\} 
\]
coincides with the root degree of $\mathscr T(\bX^{\bth}; \{U_{1}, \cdots, U_{k}\})$ for $n$ sufficiently large. 
Note that there is a finite number of vertices in $\mathscr T(\bX^{\bth}; \{U_{1}, \cdots, U_{k}\})$. Applying the previous arguments repeatedly, we see that  $\mathscr T(\sigma_{n}\bX^{n}; \{U_{1}, \cdots, U_{k}\}))$ must be the same as $\mathscr T(\bX^{\bth}; \{U_{1}, \cdots, U_{k}\})$ for $n$ sufficiently large. 
This completes the proof. 
\end{proof}

\section{Some facts about exchangeable processes}
\label{sec: a-exch}

A sequence of $n$ random variables $(\chi_{i})_{1\le i\le n}$ is said to be {\it exchangeable} if its law is unchanged by any permutation of $[n]=\{1, 2, \dots, n\}$. A process $\rY = (Y_{t})_{0\le t\le 1}\in \mathbb D([0, 1], \R)$ with $Y_{0}=0$ is {\it exchangeable} or {\it has exchangeable increments} if for all $n\in \N$, the sequence $(Y_{i/n}-Y_{(i-1)/n})_{1\le i\le n}$ is exchangeable. Kallenberg \cite{Kallenberg73} shows that any such process is necessarily of the following form:
\begin{equation}
\label{eq: canon}
Y_{t}= \alpha t+\beta b^{\br}_{t} + \sum_{i\ge 1} \theta_{i}(\mathbf 1_{\{\chi_{i}\le t\}}-t), \quad 0\le t\le 1,
\end{equation}
where $\alpha, \beta, \theta_{i}, i\ge 1$, are real-valued random variables satisfying $\beta\ge 0$, $|\theta_{1}|\ge |\theta_{2}|\ge|\theta_{3}|\ge \cdots$ and $\sum_{i}\theta_{i}^{2}<\infty$ almost surely, and are independent of the Brownian bridge $(b^{\br}_{t})_{0\le t\le 1}$ and the sequence of independent uniform variables $\chi_{i}$ on $(0, 1)$. Writing $\bth=(\theta_{i})_{i\ge 1}$, we will refer to the triple $(\alpha, \beta, \bth)$ as the {\it characteristics} of $\rY$, which is uniquely determined. Kallenberg also points out the following criterion for convergence of exchangeable processes. 

\begin{thm}[Kallenberg \cite{Kallenberg73}, Theorem 2.3]
\label{thm: cv-ex}
For each $n\in \N$, let $\rY^{n}=(Y^{n}_{t})_{t\ge 0}$ be an exchangeable process with the characteristics $(\alpha_{n}, \beta_{n}, \bth_{n})$ with $\bth_{n}=(\theta_{n, i})_{i\ge 1}$. Let $\pi_{n}=\beta_{n}\delta_{0}+\sum_{i\ge 1}\theta_{n, i}\delta_{\theta_{n, i}}$ and $\pi=\beta\delta_{0}+\sum_{i\ge 1}\theta_{i}\delta_{\theta_{i}}$. 
Then $\rY^{n}\overset{(d)}{\to} \rY$ in $\mathbb D([0, 1], \R)$ if and only if $\alpha_{n}\overset{(d)}{\to} \alpha$ and $\pi_{n}\overset{(d)}{\to} \pi$ with respect to the weak topology for finite measures on $\R$. 
\end{thm}

\begin{prop}[Knight \cite{Knight}]
\label{prop: cont-marg}
Let $\rY=(Y_{t})_{0\le t\le 1}$ be as in \eqref{eq: canon}. 
Suppose that either $\beta>0$ or $\sum_{i\ge 1}\mathbf 1_{\{\theta_{i}\ne 0\}}=\infty$. Then for all $t\in (0, 1)$, the law of $Y_{t}$ is continuous. 
\end{prop}

\begin{proof}
This is shown as an intermediate step in the proof of Lemma 1.2 \cite{Knight}. See pages 175-176 there.
\end{proof}

By replacing $Y_{t}$ with $Y_{t}-tY_{1}$, we can always bring $\alpha$ to $0$. In that case, we say that  $\rY$ is a (random) step function if the sequence $(\theta_{i})$ has at most $N\in \N$ non zero terms and $\sum_{1\le i\le N}\theta_{i}=0$ almost surely, so that 
\[
Y_{t} = \sum_{i=1}^{N} \theta_{i}\mathbf 1_{\{\chi_{i}\le t\}}, \quad 0\le t\le 1.
\]

\begin{thm}[Knight \cite{Knight}, Theorem 1.3(a) and Theorem 1.5]
\label{thm: knight}
Let $\rY=(Y_{t})_{0\le t\le 1}$ be as in \eqref{eq: canon} with $\alpha\equiv0$. Then almost surely $\rY$ has a unique infimum point if and only if $\mathbb P(\rY \text{ is a step function})=0$. 
\end{thm}

In particular, the above implies that both processes $\rY^{\bth}$ in \eqref{def: Ytheta} and $\rY^{n}$ in \eqref{def: Yn} have unique infimum points a.s. We next investigate the implication of this on the Vervaat transformation of exchangeable processes. Recall the relevant notation from around \eqref{def: Ver-op}.  

\begin{lem}
\label{lem: vervaat}
Suppose that $\rY=(Y_{t})_{0\le t\le 1}$ is an exchangeable process which has a unique infimum point $\rho_{\rY}$ a.s. Then $\rho_{\rY}$ is uniformly distributed and is independent of $\Ver(\rY)$. 
\end{lem} 

\begin{proof}
We follow the arguments below (3.14) in \cite{Bertoin01}. 
For $u\in (0, 1)$, denote by $\theta_{u}\rY$ the cyclic shift of $\rY$, namely, $\theta_{u}Y_{t}= Y_{t+u\!\!\mod 1}-Y_{u}$, $0\le t\le 1$. Since $\rY$ has a unique infimum point $\rho_{Y}$, so does $\theta_{u}\rY$, with its infimum point at $\rho_{Y}-u\!\!\mod 1$. On the other hand, we note that $\Ver( \theta_{u}(\rY))=\Ver(\rY)$. Combining this with the fact that $\theta_{u}\rY$ has the same distribution as $\rY$, we deduce that for any measurable and bounded functions $f:\R\to \R$ and $F:\mathbb D([0, 1], \R)\to \R$,
\begin{align*}
\mathbb E\Big[f\big(\rho_{Y}\big)F\big(\Ver(\rY)\big)\Big] &= \int_{0}^{1}\mathbb E\Big[f\big(\rho_{\theta_{u}Y}\big)F\big(\Ver(\theta_{u}\rY)\big)\Big] du \\
& = \mathbb E\Big[\int_{0}^{1}  f\big(\rho_{\rY}-u\big)F\big(\Ver(\rY)\big)du \Big] \\
& = \mathbb E\big[F\big(\Ver(\rY)\big)\big] \int_{0}^{1}f(u)du,
\end{align*}
where we made a change of variable in the last line. 
\end{proof}

From now on, we assume that the variables $\alpha, \beta, \bth$ are non random. Moreover, $\alpha=\beta=0, \theta_{i}\ge 0$ and $\sum_{i}\theta_{i}=\infty$. Namely, we restrict to the case of exchangeable process $\rY^{\bth}$ in \eqref{def: Ytheta}. 
Note that we have $\big(-Y^{\bth}_{(1-t)-}\big)_{0\le t\le 1}\eqd \rY^{\bth}$. In words, the law of $\rY^{\bth}$ is {\it invariant by time reversal.}
The assumption $\sum_{i}\theta_{i}=\infty$ ensures that the sample paths of $\rY^{\bth}$ has unbounded variations. In particular, the following holds true. 

\begin{prop}
\label{prop: rY-regular}
For $t\in [0, 1)$, $\mathbb P(\inf\{s>t: Y^{\bth}_{s}>Y^{\bth}_{t}\}=t)=\mathbb P(\inf\{s>t: Y^{\bth}_{s}<Y^{\bth}_{t}\}=t)=1$.  
\end{prop}

\begin{proof}
This is immediate from Theorem 1 of \cite{UB20}.
\end{proof}

\begin{prop}
\label{prop: rY-cont}
With probability 1, $\rY^{\bth}$ is continuous at its global infimum point $\rho_{\rY^{\bth}}$. 
\end{prop}

\begin{proof}
See Theorem 2 in \cite{UB20}.
\end{proof}

Let us recall that $\bX^{\bth}=\Ver(\rY^{\bth})$ and that $U_{k, 1}<U_{k, 2}< \cdots<U_{k, k}$ is the order statistics of $k$ i.i.d.~uniform variables on $(0, 1)$, which are independent of $\rY^{\bth}$.
In Section \ref{sec: cv-icrt-sp}, we have used the following properties of $\bX^{\bth}$.

\begin{lem}
\label{lem: X-prop}
We denote $m =\inf_{s\in [U_{k, 1}, U_{k, k}]}X^{\bth}_{s}$. The following events take place with probability 1:
\begin{enumerate}[(i)]
\item
$X^{\bth}_{U_{k, 1}}>0$ and $\Delta X^{\bth}_{U_{k, 1}}=0$;
\item
$0<m <X^{\bth}_{U_{k, 1}}$, so that $\tau(\bX^{\bth}, U_{k, 1}, m)<U_{k, 1}$; 
\item
$\tau(\bX^{\bth}, U_{k, 1}, m+)=\tau(\bX^{\bth}, U_{k, 1}, m)$;
\item
If $\Delta X^{\bth}(\tau(\bX^{\bth}, U_{k, 1}, m))>0$, then $X^{\bth}(\tau(\bX^{\bth}, U_{k, 1}, m)-) < m <X^{\bth}(\tau(\bX^{\bth}, U_{k, 1}, m)) $. 
\end{enumerate}
\end{lem}

\begin{proof}
\begin{enumerate}[(i)]
\item
Since $\bX^{\bth}$ only has a countable number of jumps and the law of $U_{k, 1}$ is diffuse and independent of $\bX^{\bth}$, it follows that $U_{k, 1}$ is a.s.~not a jump time of $\bX^{\bth}$. We also note that if there is some $u\in (0, 1)$ satisfying $X^{\bth}_{u-}=0$ then $\rY^{\bth}$ reaches its infimum at more than one place. By Theorem \ref{thm: knight}, the latter event has null probability. Therefore $X^{\bth}_{u}\ge X^{\bth}_{u-}>0$ for all $u\in (0, 1)$. 
\item
Suppose that $m=0$; then we can find a sequence $(t_{n})_{n\ge 1}$ contained in $[U_{k, 1}, U_{k, k}]$ and $X^{\bth}_{t_{n}}\to 0$. It follows that there exists $t_{0}\in [U_{k, 1}, U_{k, k}]$ and $X^{\bth}_{t_{0}-}=0$. By the previous arguments, this is impossible. Therefore, $m>0$. For the other inequality, let us show that for all $0< s<t<1$, 
\begin{equation}
\label{eq: inth}
\mathbb P\Big(\inf_{u\in [s, t]}X^{\bth}_{u}= X^{\bth}_{s}\Big)=0. 
\end{equation} 
We follow the arguments in the proof of Lemma 7 in \cite{Bertoin}. 
We note that on the event $\rho=\rho_{\rY^{\bth}}<1- t$, the interval $[s, t]$ is shifted to $[\rho+s, \rho+t]$ in $\rY^{\bth}$. Combined with Lemma~\ref{lem: vervaat}, we deduce that 
\begin{align*}
(1-t)\,\mathbb P\Big(\inf_{u\in [s,t]}X^{\bth}_{u}= X^{\bth}_{s}\Big)  & =\mathbb E\int_{0}^{1-t} \mathbf 1_{\{\inf_{u\in [s,t]}X^{\bth}_{u}= X^{\bth}_{s}; \rho=v\}} dv \\
& =\mathbb E\int_{0}^{1-t}\mathbf 1_{\{\inf_{u\in [v+s,v+t]}Y^{\bth}_{u}= Y^{\bth}_{v+s}; \rho=v\}} dv \\
& \le \int_{0}^{1-t} \mathbb P\Big(\inf_{u\in [v+s, v+t]}Y^{\bth}_{u}=Y^{\bth}_{v+s}\Big) dv\\
& =  \int_{0}^{1-t} \mathbb P\Big(\inf_{u\in [s, t]}Y^{\bth}_{u}=Y^{\bth}_{s}\Big) dv = 0,
\end{align*}
where we have relied on $s>0$ in the penultimate line, then used exchangeability and Lemma~\ref{prop: rY-regular} for the last line. This proves \eqref{eq: inth} and the desired result follows. 

\item
Again, it suffices to prove the statement for fixed $0<s<t<1$. 
We observe that on the event $\tau':=\tau(\bX^{\bth}, s, m_{s, t}+) >\tau:= \tau(\bX^{\bth}, s, m_{s, t})\le s$, we will have $\inf_{u\in [\tau'-\epsilon,\tau']}X^{\bth}_{u}=m$ for all $0<\epsilon<q'-q$. It follows that we can find some rationals $q<q'$ so that $\bX^{\bth}$ restricted to $[q, q']$ attains minimum at two different locations. 
Arguing as previously, we see that this implies $Y^{\bth}$ restricted to some interval $[r, r']$ will attain minimum at two different locations, with $0<r<r'<1$. But $Y^{\bth}$ restricted to $[r, r']$ is still an exchangeable process and has a similar representation as in \eqref{eq: canon}. It follows from Theorem \ref{thm: knight} this event has null probability. 

\item
As before, fix $0<s<t<1$ and denote $m_{s, t}=\inf_{u\in [s, t]}X^{\bth}_{u}$, $\tau = \tau(\bX^{\bth}, s, m_{s, t})$. By definition, $X^{\bth}_{\tau-}\le m_{s, t}$ and $X^{\bth}_{\tau}\ge m_{s, t}$; so we only need to exclude the possibilities that $X^{\bth}_{\tau-}=m_{s, t}$ or $X^{\bth}_{\tau}=m_{s, t}$. Note the only jump times of $\rY^{\bth}$ are $\chi_{i}$, $i\ge 1$, which are independently and uniformly distributed. Let us introduce
\[
\widetilde Y^{\bth, i}_{t}= \sum_{j\ne i} \theta_{j}(\mathbf 1_{\{\chi_{j}\le t\}}-t) = Y^{\bth}_{t}-\theta_{i}(\mathbf 1_{\{\chi_{i}\le t\}}-t), \quad 0\le t\le 1.
\]
Note that $\chi_{i}$ is independent of $(\widetilde Y^{\bth, i}_{t})_{0\le t\le 1}$ and we have $\widetilde Y^{\bth, i}_{\chi_{i}}=Y^{\bth}_{\chi_{i}-}+\theta_{i}\chi_{i}$, as well as $\widetilde Y^{\bth, i}_{\chi_{i}}=Y^{\bth}_{\chi_{i}-}-\theta_{i}(1-\chi_{i})$. Thanks to Proposition \ref{prop: cont-marg}, we have for any $z\in \R$, 
\[
\mathbb P\big(\exists\, u\le s: Y^{\bth}_{u-}=z<Y^{\bth}_{u}\big) \le 
\sum_{i\ge 1}\mathbb P\big(\widetilde Y^{\bth, i}_{\chi_{i}}=z+\theta_{i}\chi_{i} \big) = \sum_{i\ge 1}\int_{0}^{1}\mathbb P\big(\widetilde Y^{\bth, i}_{u}=z+u \big) du =0.  
\]
Similarly, 
\[
\mathbb P\big(\exists\, u\le s: Y^{\bth}_{u-}<Y^{\bth}_{u}=z\big)\le \sum_{i\ge 1}\int_{0}^{1}\mathbb P\big(\widetilde Y^{\bth, i}_{u}=z-\theta_i u \big) du = 0. 
\]
Let us note that the above also holds when conditioned on $Y^{\bth}_{s}$, since $(Y^{\bth}_{us})_{0\le u\le 1}$ is an exchangeable process. Moreover, conditioning on $Y^{\bth}_{s}$ and the respective subsets of $\chi_{i}$'s that are contained in $[0, s]$ and $[s, 1]$,  $(Y^{\bth}_{u})_{u\in [0, s]}$ and $(Y^{\bth}_{u})_{u\in [s, 1]}$ are independent, as a consequence of \eqref{eq: canon}. 
It follows that $\inf_{v\in [s, t]}Y^{\bth}_{v}$ is conditionally independent of $(Y^{\bth}_{u})_{u\in [0, s]}$. 
We find via integration that
\[
\mathbb P\Big(\exists\, u\le s: Y^{\bth}_{u-}=\inf_{v\in [s, t]}Y^{\bth}_{v} <Y^{\bth}_{u}\Big) =\mathbb P\Big(\exists\, u\le s: Y^{\bth}_{u-}<Y^{\bth}_{u}=\inf_{v\in [s, t]}Y^{\bth}_{v} \Big)= 0. 
\]
Arguing as previously, we conclude this holds similarly for $\bX^{\bth}$. 
\end{enumerate}
\end{proof}

Recall $\mathrm Z=(Z_{t})_{0\le t\le 1}$ from \eqref{def: Z} and $\tilde U_{k, i}$ from \eqref{def: Utilde}. Let us denote $g_i=g(\mathrm Z, \tilde U_{k, i})$ and $d_i=d(\mathrm Z, \tilde U_{k, i})$, $1\le i\le k$. 
\begin{lem}
\label{lem: Z-prop}
For each $1\le i\le k$, with probability $1$, we have $g_i<\tilde U_{k, i}<d_i$ and $\inf_{u\in [0, g_i-\epsilon]}Z_u >\inf_{u\in [0, \tilde U_{k, i}]}Z_{u}>\inf_{u\in [0, d_i+\epsilon]}Z_u$ for all $\epsilon >0$. 
\end{lem}

\begin{proof}
Using arguments similar to the ones leading to \eqref{eq: inth} and combining them with the time reversal property of $\rY^{\bth}$, we can show that for all $0\le s<t<1$, 
\begin{equation}
\label{eq: ineqe}
\mathbb P\Big(\inf_{u\in [s, t]}X^{\bth}_{u}= X^{\bth}_{t-}\Big)=0. 
\end{equation}
If $g_i=\tilde U_{k, i}$, then $Z_{\tilde U_{k, i}-}=\inf_{u\in [0,\tilde U_{k, i}]}Z_u$, from which it follows $X^{\bth}_{U_{k, i}-}=\inf_{u\in [b, U_{k, i}]}X^{\bth}_u$. 
This is clearly impossible because of \eqref{eq: ineqe} and the fact that the law of $U_{k, i}$ is independent of $\bX^{\bth}$. 
Similarly, we can argue that $\tilde U_{k, i}<d_i$ a.s. 
If there is some $\epsilon>0$ satisfying $\inf_{u\in [0, g_i-\epsilon]}Z_u =\inf_{u\in [0, U_{k, i}]}Z_{u}$, then $\bX^{\bth}$ attains a local minimum at two different locations. We have seen in the proof of Lemma \ref{lem: X-prop} this occurs with null probability.  
\end{proof}

\appendix
\section{Convergence of vertex degrees to local times in a stable tree}
\label{sec: A1}

Recall the sequence of spanning trees $\cT_{k}$ of the $\alpha$-stable tree $\cT_{\alpha}$. 
We prove here the approximation  \eqref{def: loc-stable''} for the local times of the branch points. 

\begin{prop}
\label{prop: a-loc}
With probability $1$, we have that
\begin{equation}
\label{def: loc-stable3}
\Delta^{\alpha}(b)=\lim_{k\to\infty}\frac{\deg(b, \cT_{k})}{k^{1/\alpha}} \quad\text{a.s.}
\end{equation}
holds for all $b\in \Br(\cT_{\alpha})$. 
\end{prop}



\paragraph{Poissonian marking.}
Our approach here makes use of the Poissonian marking from \cite{DuLG02}.  
Recall the canonical process $\mathrm e=(e_{s})_{s\ge 0}$ of the Skorokhod space and its lifetime $\zeta=\zeta(e) =\inf\{t>0: e_{s}= 0\,\forall\,s\ge t\}\in (0, \infty)$. Let $U_{1}<U_{2}<U_{3}<\cdots < U_{N(\lambda)}$ be the jumps of a Poisson process on $[0, \zeta]$ of rate $\lambda>0$ per unit time. Standard properties of Poisson processes imply that
\begin{itemize}
\item
$N(\lambda)$ has the Poisson distribution of mean $\lambda \zeta$.
\item
Given $N(\lambda)=k$, $k\in\N$, $\zeta^{-1}(U_{1}, U_{2}, \dots, U_{N(\lambda)})$ is distributed as the order statistics of $k$ independent uniform variables on $(0, 1)$. 
\end{itemize}
Write $\mathbf U=(U_{1}, U_{2}, \cdots, U_{N(\lambda)})$. 
We define $D(\mathrm e\,; \mathbf U)$ to be the number of distinct values in the following collection:
\[
\inf_{0\le s\le u} e_{s}, \quad u \in \mathbf U. 
\]
In other words, $D(\mathrm e\,; \mathbf U)$ simply counts the number of excursion intervals of $\mathrm e$ above its infimum which contain at least a mark from the Poisson process. 
For $x>0$, let $\mathbb P_{x}$ denote the law of $\rY^{\alpha}$ stopped when it first reaches the level $-x$. The excursion theory implies that under $\mathbb P_{x}$, $D(\mathrm e\,; \mathbf U)$ is distributed as a Poisson variable of mean
\[
x\cdot\mathbb N(1-e^{-\lambda \zeta}) = x\lambda^{1/\alpha},
\]
where we have used \eqref{eq: zeta-dis}. 
It follows that 
\begin{equation}
\label{lim: deg1}
\lim_{\lambda\to\infty} \frac{D(\mathrm e\,; \mathbf U)}{\lambda^{1/\alpha}} = x, 
\end{equation}
where the limit holds $\mathbb P_{x}$-almost surely.

\paragraph{Large jumps in an excursion process. } Our second ingredient for the proof of Proposition \ref{prop: a-loc} is a description of $\mathrm e$ under $\mathbb N$ conditioned to have at least $k$ jumps of size at least $\epsilon$. 
Recall that the jumps of the stable process $\mathrm Y^{\alpha}$ follow a Poisson point process of intensity measure $\pi(dx) = c_{\alpha}x^{-\alpha-1}dx\mathbf{1}_{\{x>0\}}$, with $c_{\alpha}= \alpha(\alpha-1)/\Gamma(2-\alpha)$. 
In particular, the first moment that $\mathrm Y^{\alpha}$ has a jump that is at least $\epsilon$ large, namely, 
\[
\gamma_{1}(\mathrm Y^{\alpha})=\inf\{s>0: \Delta Y^{\alpha}_{s}\ge \epsilon\}, 
\]
is distributed as an exponential variable of rate $\pi([\epsilon, \infty))$. Recall from \eqref{def: sigma-time} the notation $\sigma_{\rY^{\alpha}}(\gamma_{1}(\rY^{\alpha}))$. The strong Markov property implies that conditional on $\Delta Y^{\alpha}_{\gamma_{1}(\mathrm Y^{\alpha})}=x$, the process $\mathrm Z^{1}$ defined by
\[
Z^{1}_{s}=Y^{\alpha}_{s+\gamma_{1}(\mathrm Y^{\alpha})}, \quad 0\le s \le \sigma_{\rY^{\alpha}}(\gamma_{1}(\rY^{\alpha}))-\gamma_{1}(\mathrm Y^{\alpha}),
\]
has the law $\mathbb P_{x}$. Iterate this procedure, we can obtain the same description for each of the portions of $\mathrm Y^{\alpha}$ between its $n$-th moment of having a jump at least $\epsilon$ large and the moment the process falls back to the level prior to this jump. Combining this with the excursion theory, we find the following

\begin{lem}
\label{lem: a-decomp}
Fix $n\in \N$ and $\epsilon\in (0, \infty)$. Set $\gamma_{0}=0$ and for $1\le i\le n$, let
\[
\gamma_{i}=\inf\{s>\gamma_{i-1}: \Delta e_{s}\ge \epsilon\}, \quad \text{and if }\gamma_{i}<\infty, \text{ set } \sigma_{i}=\inf\{s>\gamma_{i}: e_{s}\le e_{\gamma_{i}-}\}.
\]
If $\gamma_{i}<\infty$, define also the process $\mathrm Z^{i}=(Z^{i}_{s})$ by 
\[
Z^{i}_{s}=e_{s+\gamma_{i}}, \quad 0\le s \le \sigma_{i}-\gamma_{i}.
\]
Then for $1\le i\le n$, under $\mathbb N(\cdot\,|\, \gamma_{i}<\infty)$, $\mathrm Z^{i}$ has the law $\mathbb P_{\Delta e_{\gamma_{i}}}$. 
\end{lem}

\begin{proof}[Proof of Proposition \ref{prop: a-loc}]
Recall the real tree $\mathcal T_{\alpha}$ encoded by the height process $\mathrm H$, which is itself a function of the excursion $\mathrm e$ under $\mathbb N$. 
Let $\mathcal P(\lambda)=(J_{i})_{1\le i\le L(\lambda)}$ be a Poisson process on $(0, \zeta)$ with rate $\lambda$ per unit time. Denote by $\cT(\lambda)$ the subtree of $\cT$ spanned by $J_{1}, J_{2}, \cdots, J_{L(\lambda)}$. Clearly, on the event that $L(\lambda)=k$, $\cT(\lambda)$ has the same distribution as $\cT_{k}$ under $\mathbb N$.  Suppose that $\mathrm e$ has at least $n$ jumps that are at leat $\epsilon$ large, the first $n$ of which occur respectively at $\gamma_{1}, \gamma_{2}, \cdots, \gamma_{n}$, and are of respective sizes $\Delta_{1}, \Delta_{2}, \cdots, \Delta_{n}$.  For $1\le i\le n$, let
\[
\sigma_{i}=\inf\{s>\gamma_{i}: e_{s}\le e_{\gamma_{i}-}\} \quad \text{and}\quad Z^{i}_{s} = e_{s+\gamma_{i}}, 0\le s\le \sigma_{i}-\gamma_{i}.
\]
Then Lemma \ref{lem: a-decomp} says that each $\mathrm Z^{i}=(Z^{i}_{s})_{0\le s\le \sigma_{i}-\gamma_{i}}$ has the law $\mathbb P_{\Delta_{i}}$. On the other hand, the jump at $\gamma_{i}$ corresponds to a unique branch point $p(\gamma_{i})$ of $\cT$. Let's find its degree in the reduced tree $\cT(\lambda)$. 
Let $I_{i}=\{1\le i\le L(\lambda): \gamma_{i}<J_{i}<\sigma_{i}\}$, the set of indices for the Poisson marks that fall into $(\gamma_{i}, \sigma_{i})$. By the properties of Poisson processes, the subset of marks $\{J_{i}: \gamma_{i}<J_{i}<\sigma_{i}\}$ has the distribution of a Poisson point process of intensity $\lambda$ on $(\gamma_{i}, \sigma_{i})$. 
From the definition of the reduced tree, we see that 
\[
\Big(\deg\big(p(\gamma_{i}), \cT(\lambda)\big)-1\Big)_{+} = D\big(\mathrm X^{i}; (J_{i}: i\in I_{i})\big)
\]
Then Lemma \ref{lem: a-decomp} together with \eqref{lim: deg1} shows that 
\[
\lim_{\lambda\to\infty} \frac{\deg\big(p(\gamma_{i}), \cT(\lambda)\big)}{\lambda^{1/\alpha}} = \Delta_{i}, \quad 1\le i\le n,
\]
almost surely under $\mathbb N(\cdot \,|\, \gamma_{i}<\infty)$. Since $n\in \N$ and $\epsilon>0$ are arbitrary, this allows us to conclude that for all $s$ such that $\Delta e_{s}>0$, we have
\[
\lim_{\lambda\to\infty} \frac{\deg\big(p(s), \cT(\lambda)\big)}{\lambda^{1/\alpha}} = \Delta e_{s}, \quad \mathbb N\text{-a.e.}
\]
Thanks to \eqref{eq: normalise}, we find that the above limit also holds $\nr$-a.s. On the other hand, since $L(\lambda)$ has the Poisson distribution with mean $\zeta\lambda$, we have $L(\lambda)/\lambda\to \zeta$ as $\lambda\to\infty$, $\mathbb N$-a.e. 
Applying once again \eqref{eq: normalise}, we obtain that 
\[
\lim_{\lambda\to\infty} \frac{L(\lambda)}{\lambda} = 1, \quad \nr\text{-a.s.}
\]
Combined with the previous limit, this implies that for all $s$ such that $\Delta e_{s}>0$, we have
\[
\lim_{\lambda\to\infty} \frac{\deg\big(p(s), \cT(\lambda)\big)}{L(\lambda)^{1/\alpha}} = \Delta e_{s}, \quad \nr\text{-a.s.}
\]
In particular, the above limit holds along the (random) subsequence $(\lambda_{k})_{k\ge 1}$, where $\lambda_{k}=\min\{\lambda: L(\lambda)\ge k\}$. The desired result follows. 
\end{proof}


\section{Some facts about Skorokhod's topology}
\label{sec: topo}

Here, we gather some results on Skorokhod's topology used in the proof of Proposition~\ref{prop: cv-sp}. 
We denote by $\mathbb D=\mathbb D([0, 1], \R)$ the space of c\`adl\`ag functions defined on $[0, 1]$ equipped with Skorokhod's topology. 
In the sequel, $\|\cdot\|$ stands for the uniform norm on $[0, 1]$, and $\mathrm{Id}$ the identity map of $[0, 1]$. 

\begin{lem}
\label{lem: D15}
Suppose that $\mathbf x_{n}\to \mathbf x$ in $\mathbb D$ and $b_{n}\to b\in (0, 1)$. Suppose that either $\mathbf x$ is continuous at $b$, or $\Delta x_{n}(b_{n})\to \Delta x(b)\ne 0$. Then 
\[
\Big\{x_{n}(b_{n}+(1-b_{n})t)-x_{n}(b_{n}): t\in [0, 1]\Big\} \xrightarrow{n\to\infty} \Big\{x(b+(1-b)t)-x(b): t\in [0, 1]\Big\} \quad \text{ in $\mathbb D$}.
\]
\end{lem}

\begin{proof}
In the first place, let us assume that $\mathbf x$ is continuous at $b$. Since $\mathbf x_{n}\to \mathbf x$ in $\mathbb D$, we can find a sequence of strictly increasing and continuous bijections $\lambda_{n}: [0, 1]\to [0, 1]$ so that $\|\lambda_{n}-\mathrm{Id}\|\to 0$ and $\|\mathbf x_{n}-\mathbf x\circ \lambda_{n}\|\to 0$. Meanwhile, $b$ being a point of continuity for $\mathbf x$, there is a sequence of positive real numbers  $\delta_{n}\to 0$ so that 
\[
\sup_{s\in [b-\delta_{n}, b+\delta_{n}]}|x(s)-x(b)|\le \frac1n. 
\]
Denote by $q_{n}$ the real number satisfying $(1-b)q_{n}=\lambda_{n}(b_{n}+(1-b_{n})\delta_{n})-b$. The fact that $b_{n}\to b$ implies that $q_{n}-\delta_{n}\to 0$. Therefore, it is possible to find, at least for $n$ sufficiently large, a strictly increasing and continuous bijection $f_{n}: [0, \delta_{n}]\to [0, q_{n}]$ satisfying $\sup_{u\in [0, \delta_{n}]}|f_{n}(u)-u|\to 0$. We define the function $\tilde\lambda_{n}:[0, 1]\to [0, 1]$ as follows: if $u\le \delta_{n}$, $\tilde\lambda_{n}(u)=f_{n}(u)$; if $\delta_{n}<u\le 1$, let $\tilde\lambda_{n}(u)$ be defined by
\[
b+ (1-b)\tilde\lambda_{n}(u)=\lambda_{n}\big(b_{n}+(1-b_{n})u\big).
\]
It can be readily checked that $\tilde\lambda_{n}:[0, 1]\to [0, 1]$ is strictly increasing, bijective, continuous, and satisfies $\|\tilde\lambda_{n}-\mathrm{Id}\|\to 0$. Moreover, we have
\begin{align*}
&\sup_{t\in [0, 1]}\Big|x_{n}\big(b_{n}+(1-b_{n})t\big)-x_{n}(b_{n})-x\big(b+(1-b)\tilde\lambda_{n}(t)\big)+x(b)\Big|\\
&\qquad\qquad\qquad\le\sup_{t\in [0, 1]} \Big|x_{n}\big(b_{n}+(1-b_{n})t\big)-x\circ \lambda_{n}\big(b_{n}+(1-b_{n})t\big)\Big| + |x_{n}(b_{n})-x(b)|\\
&\qquad\qquad\qquad\qquad+\sup_{t\in [0, 1]}\Big|x\circ\lambda_{n}\big(b_{n}+(1-b_{n})t\big)-x\big(b+(1-b)\tilde\lambda_{n}(t)\big)\Big|,
\end{align*}
where the first term on the second line is bounded by $\|\mathbf x_{n}-\mathbf x\circ\lambda_{n}\|$, the second term tends to $0$ as $\mathbf x$ is continuous at $b$, and the term on the last line is $\le 1/n$ for $n$ sufficiently large, by the choice of $\tilde\lambda_{n}$. This proves the statement when $\mathbf x$ is continuous at $b$. If, instead, $\Delta x(b)\ne 0$, then define $\tilde{x}(t) = x(t)-\Delta x(b)\mathbf 1_{\{t\ge b\}}$, $t\in [0, 1]$, so that $\Delta\tilde x(b)=0$. Similarly, let $\tilde{x}_{n}(t) = x_{n}(t)-\Delta x_{n}(b_{n})\mathbf 1_{\{t\ge b_{n}\}}$, $t\in [0, 1]$. We can show that $(\tilde x_{n}(t))_{t\in [0, 1]}\to (\tilde x(t))_{t\in [0, 1]}$ in $\mathbb D$. The conclusion follows, as $x(u) = \tilde x(u)+\Delta x(b)$ for all $u\ge b$, and similarly for $\mathbf x_{n}$. 
\end{proof}

Let us recall the notation $m(\mathbf x, s, t)=\inf_{u\in [s, t]}x(u)$. A proof of the following lemma can be found for instance in Lemma B.3 of \cite{BrDuWa21}. 

\begin{lem}
\label{lem: D-1}
Suppose that $\mathbf x_n\to \mathbf x$ in $\mathbb D$ and that $0\le s<t\le 1$ satisfy $\Delta x(s)=\Delta x(t)=0$. Then for any $s_n\to s$ and $t_n\to t$, we have $m(\mathbf x_n, s_n, t_n) \to m(\mathbf x, s, t)$. 
\end{lem}

\begin{lem}
\label{lem: D0}
For each $n\in \N$, let $\mathbf y_{n}\in \mathbb D$ satisfy $y_{n}(1)=y_{n}(0)=0$. Let $\Ver(\mathbf y_{n})$ be as defined in \eqref{def: Ver-op}. Suppose that $\mathbf y\in \mathbb D$ with $y(1)=y(0)=0$ has a unique infimum point $\rho_{\mathbf y}$ and $\mathbf y$ is continuous at $\rho_{\mathbf y}$. Then $\mathbf y_{n}\to \mathbf y$ in $ \mathbb D$ implies $\Ver(\mathbf y_{n})\to \Ver(\mathbf y)$ in $ \mathbb D$. 
\end{lem}

\begin{proof}
Let us first show that $\rho_{\mathbf y_n}=\inf\{t>0: y_n(t-)\wedge y(t)=\inf_{s\in [0, 1]}y_n(s)\}\to \rho_{\mathbf y}$, as $n\to\infty$. Indeed, since $\mathbf y$ has a unique infimum point and since $y(0)=y(1)=0$, we must have $\rho_y\in (0, 1)$ and $m(\mathbf y, 0, 1)<0$. Moreover, for all $\epsilon>0$, there exists some $t_\epsilon\in (\rho_{\mathbf y}, \rho_{\mathbf y}+\epsilon)$ so that $\Delta y(t_\epsilon)=0$ and $0>m(\mathbf y, t_\epsilon, 1)>m(\mathbf y, 0, 1)$. Meanwhile, as $\mathbf y_n\to\mathbf y$ and $t_\epsilon$ is a continuity point of $\mathbf y$, we have $m(\mathbf y_n, t_\epsilon, 1)\to m(\mathbf y, t_\epsilon, 1)$ and $m(\mathbf y_n, 0, 1)\to m(\mathbf y, 0, 1)$, so that $m(\mathbf y_n, t_\epsilon, 1) > m(\mathbf y, 0, 1)$ and consequently $\rho_{\mathbf y_n}\le t_\epsilon\le \rho_{\mathbf y}+\epsilon$ for $n$ large enough. Similarly we can show that $\limsup_{n\to\infty}\rho_{\mathbf y_n}\ge \rho_{\mathbf y}-\epsilon$ for any $\epsilon>0$. Combining the two arguments, we deduce that $\rho_{\mathbf y_n}\to \rho_{\mathbf y}$. The convergence of $\Ver(\mathbf y_{n})$ to $\Ver(\mathbf y)$ can be shown by arguments similar to the proof of Lemma \ref{lem: D15}. 
\end{proof}

For $\mathbf x=(x(t))_{t\in [0, 1]}\in \mathbb D$, $t\ge 0$ and $r\in \R$, recall that $\tau(\mathbf x, t, r)=\inf\{s\le t: \inf_{u\in [s, t]}x(u)\ge r\}$, with the convention that $\inf\varnothing=\infty$. Note that $r\mapsto \tau(\mathbf x, t, r)$ is non decreasing and left-continuous. 

\begin{lem}
\label{lem: D1}
Suppose that $\mathbf x_{n}\to \mathbf x$ in $\mathbb D$. 
Assume that $\tau(\mathbf x, t, r+)=\tau(\mathbf x, t, r)<t$ and $\Delta x(t)=0$. Then for all $(t_n, r_{n})\to (t, r)$, we have $\tau(\mathbf x_{n}, t_{n}, r_{n})\to \tau(\mathbf x, t, r)$. 
In the case where $\Delta x(\tau(\mathbf x, t, r))>0$, assume further $x(\tau(\mathbf x, t, r)-) <r<x(\tau(\mathbf x, t, r))$; 
 then we also have $\Delta x_{n}(\tau(\mathbf x_{n}, t_{n}, r_{n}))\to \Delta x(\tau(\mathbf x, t, r))$. 
\end{lem} 

\begin{proof}
 Let us denote $\tau_0=\tau(\mathbf x, t, r)$. By definition, for all $\epsilon>0$, there exists some $t_\epsilon\in (\tau_0-\epsilon, \tau_0)$ satisfying $\Delta x(t_\epsilon)=0$ and $x(t_\epsilon)<r$. This implies $x_n(t_\epsilon)<r_n$ for $n$ large enough, and subsequently  $\tau(\mathbf x_n, t_n, r_n)\ge t_\epsilon\ge \tau_0-\epsilon$. Hence, $\liminf_{n\to\infty}\tau(\mathbf x_n, t_n, r_n)\ge \tau_0$. On the other hand, as $\tau(\mathbf x, t, r+)=\tau(\mathbf x, t, r)$, for all $\epsilon>0$, we can find some $r'>r$ so that $\tau(\mathbf x, t, r')<\tau_0+\epsilon$. This implies that there exists some $s_\epsilon\in (\tau_0, \tau_0+\epsilon)$ satisfying $\Delta x(s_\epsilon)=0$ and $m(\mathbf x, s_\epsilon, t)\ge r'>r$. It follows that $m(\mathbf x_n, s_\epsilon, t_n)\ge r_n$ for $n$ large enough. Hence, $\tau(\mathbf x_n, t_n, r_n)\le s_\epsilon\le \tau_0+\epsilon$. This shows that $\tau(\mathbf x_{n}, t_{n}, r_{n})\to \tau(\mathbf x, t, r)$. 

Suppose from now on that $\Delta x(\tau_0)>0$ and $x(\tau_{0}-)<r<x(\tau_{0})$. Since $\mathbf x_n\to \mathbf x$, there is a  sequence of strictly increasing and continuous bijections $\lambda_{n}: [0, 1]\to [0, 1]$ so that $\|\lambda_{n}-\mathrm{Id}\|\to 0$ and $\|\mathbf x_{n}-\mathbf x\circ \lambda_{n}\|\to 0$. Denote by $\tau_{n}=\lambda_{n}^{-1}(\tau_{0})$. Note that we have $ x_{n}(\tau_{n})\to  x(\tau_{0})$ and $x_{n}(\tau_{n}-)\to x(\tau_{0}-)$. The statement will follow once we show that 
\begin{equation}
\label{eq: taun}
\tau(x_{n}, t_{n}, r_{n})=\tau_{n} \quad\text{ for $n$ sufficiently large. }
\end{equation}
{\it Proof of \eqref{eq: taun}. } 
Since $x_{n}(\tau_{n}-)\to x(\tau_{0}-)<r$, we have $x_{n}(\tau_{n}-)<r_{n}$ for $n$ sufficiently large. Hence, $\tau(\mathbf x_{n}, t_{n}, r_{n})\ge \tau_{n}$. 
Next, let $\epsilon>0$. We can find some $\delta>0$ so that $\sup_{\tau_{0}\le u\le \tau_{0}+\delta} |x(u)-x(\tau_{0})|<\epsilon/2$. 
Combining this with  $\|\mathbf x_{n}\circ\lambda^{-1}_{n}-\mathbf x\|\to 0$, we find that 
\[
\sup_{\tau_{n}\le u\le \tau_{n}+\delta/2}|x_{n}(u)-x(\tau_{0})|\le \epsilon \ \text{ for $n$ sufficiently large.}
\]
It follows that $m(\mathbf x_{n}, \tau_{n}, \tau_{n}+\delta/2) \ge x(\tau_{0})-\epsilon$. Since $x(\tau_{0})>r$ and $r_{n}\to r$, we deduce that  $m(\mathbf x_{n}, \tau_{n}, \tau_{n}+\delta/2) \ge
r_{n}$ for $n$ sufficiently large. Moreover, from the previous part of the proof for the lemma, we have $\liminf_{n\to\infty}m(\mathbf x_{n}, \tau_{n}+\delta/2, t_{n})\ge m(\mathbf x, \tau_{0}+\delta/4, t)>r$. Putting the two together, we have $m(\mathbf x_{n}, \tau_{n}, t_{n})\ge r_{n}$ for $n$ sufficiently large, from which it follows that $\tau(\mathbf x_{n}, t_{n}, r_{n})\le \tau_{n}$. This completes the proof of \eqref{eq: taun}. 
\end{proof}

For $\mathbf x=(x(t))_{t\in [0, 1]}\in \mathbb D$, let us denote $\underline x(t)=\inf_{s\in[0, t]}x(s)$ and $\underline{\mathbf x}=(\underline x(t))_{t\in [0, 1]}\in\mathbb D$.
For $t\in (0, 1)$, let us define $g(\mathbf x, t)$ and $d(\mathbf x, t)$ as follows: 
\[
g(\mathbf x, t) = \sup\{ s\le t: x(s)\wedge x(s-) =  \underline x(t)\} \quad \text{and} \quad d(\mathbf x, t)=\inf\{s>t: x(s) \wedge x(s-) \le \underline x(t)\}\wedge 1.
\]
\begin{lem}
\label{lem: D2}
Suppose that $\mathbf x_{n}\to \mathbf x$ in $\mathbb D$ and $t_{n}\to t\in (0, 1)$. Suppose further that $t\mapsto \underline x(t)$ is continuous and  $g(\mathbf x, t)<t<d(\mathbf x, t)$. Moreover, for every $\epsilon>0$, we have $\underline x(g(\mathbf x, t)-\epsilon)>\underline x(t)>\underline x(d(\mathbf x, t)+\epsilon)$. Then 
$g(\mathbf x_{n}, t_{n})\to g(\mathbf x, t)$ and $d(\mathbf x_{n}, t_{n})\to d(\mathbf x, t)$.
\end{lem}

\begin{proof}
We will use the shorthand notation $g:=g(\mathbf x, t)$ and $d:=d(\mathbf x, t)$. 
Since $\underline{\mathbf x}$ is continuous, we deduce from $\mathbf x_{n}\to \mathbf x$ that $\underline x_{n}(s_{n})\to \underline x(s)$ for any $s_{n}\to s$. In consequence, $\underline x_{n}(g-\epsilon) \to \underline x(g-\epsilon) > \underline x(t) = \lim \underline x_{n}(t_{n})$. It follows that $g-\epsilon\le g(\mathbf x_{n}, t)$. Since this is true for any $\epsilon$, we deduce that $\liminf_{n\to\infty}g(\mathbf x_{n}, t)\ge g$. Meanwhile, if $s\in (g, g+\epsilon)$ is a point of continuity for $\mathbf x$, then we must have $x_{n}(s)\to x(s)>\underline x(t)$. Hence, $\limsup_{n\to\infty} g(\mathbf x_{n}, t)\le g+\epsilon$ for all $\epsilon>0$. This shows that $g(\mathbf x_{n}, t_{n})\to g$. The proof for $d(\mathbf x_{n}, t_{n})$ is similar. 
\end{proof}

{\small
\setlength{\bibsep}{.2em}

}

\end{document}